\documentclass[11pt,a4paper, final, twoside]{article}
\usepackage[style=numeric-comp, backend=bibtex, doi=false, isbn=false, natbib=true]{biblatex}
\usepackage{filecontents}
\usepackage[centertags]{amsmath}
\usepackage{amsfonts}
\usepackage{amssymb}
\usepackage{eqparbox}
\usepackage{amsthm}
\usepackage{newlfont}
\usepackage{fancyhdr}
\usepackage{amscd}
\usepackage{scalerel}
\usepackage{graphicx}
\usepackage{subfigure}
\usepackage{afterpage}
\usepackage{relsize}
\usepackage{sectsty}
\usepackage{siunitx}
\usepackage{etoolbox}
\usepackage{caption}
\usepackage{array}
\usepackage{booktabs}
\usepackage[table]{xcolor}
\usepackage{xcolor,colortbl}
\usepackage{colortbl,multirow,hhline}
\definecolor{blue(ncs)}{rgb}{0.0, 0.53, 0.74}
\usepackage[colorlinks=true, urlcolor=blue(ncs), linkcolor=blue(ncs), citecolor=blue(ncs)]{hyperref}
\usepackage{euscript,pifont,mathrsfs,latexsym,graphicx}
\usepackage{mathrsfs}
\usepackage{latexsym,epstopdf,float}
\usepackage[atend]{bookmark}
\bookmarksetup{
  open,
  openlevel=2,
  numbered,
  addtohook={%
    \ifnum\bookmarkget{level}>1 %
      \DisableBookmarkNumbering
    \fi
  },
}
\makeatletter
\newcommand*{\DisableBookmarkNumbering}{%
  \let\numberline\@gobble
}
\makeatother
\patchcmd\thebibliography{\labelsep}{\labelsep\itemsep=0pt\parsep=0pt\relax}{}{\typeout{Couldn't patch the command}}

\theoremstyle{plain}
\newtheorem{thm}{Theorem}[section]
\newtheorem{cor}[thm]{Corollary}
\newtheorem{lem}[thm]{Lemma}

\theoremstyle{definition}
\newtheorem{defn}{Definition}[section]
\newtheorem*{assum}{Assumptions}
\theoremstyle{remark}
\newtheorem{rem}{Remark}[section]
\newtheorem{exa}[rem]{Example}

 \numberwithin{equation}{section}

 \newcommand*{\medcap}{\mathbin{\scalebox{1.2}{\ensuremath{\bigcap}}}}
 \newcommand*{\medcup}{\mathbin{\scalebox{1.3}{\ensuremath{\bigcup}}}}
 \newcommand*{\sOmega}{\mathbin{\scalebox{.5}{\ensuremath{\Omega}}}}

\begin{document}

\title{Approximate solutions for robust multiobjective optimization programming in Asplund spaces}
\author{Maryam Saadati$^*$, Morteza Oveisiha}
\date{}

\pdfbookmark[0]{Approximate solutions for robust multiobjective optimization programming in Asplund spaces}{frontmatter}
\maketitle
\begin{center}
Department of Pure Mathematics, Faculty of Science, Imam Khomeini International University, P.O. Box 34149-16818, Qazvin, Iran.
\\
E-mail: m.saadati@edu.ikiu.ac.ir, oveisiha@sci.ikiu.ac.ir
\end{center}

\pdfbookmark[section]{\abstractname}{abstract}
\begin{abstract}
In this paper, we study a nonsmooth/nonconvex multiobjective optimization problem with uncertain constraints in arbitrary Asplund spaces. We first provide necessary optimality condition in a fuzzy form for approximate weakly robust efficient solutions and then establish necessary optimality theorem for approximate weakly robust quasi-efficient solutions of the problem in the sense of the limiting subdifferential by exploiting a fuzzy optimality condition in terms of the Fr\'{e}chet subdifferential. Sufficient conditions for approximate (weakly) robust quasi-efficient solutions to such a problem are also driven under the new concept of generalized pseudo convex functions. Finally, we address an approximate Mond-Weir-type dual robust problem to the reference problem and explore weak, strong, and converse duality properties under assumptions of pseudo convexity.
\end{abstract}
\textbf{Keywords}\hspace{3mm}Approximate solutions . Optimality conditions . Duality . Limiting subdifferential . Generalized convexity . Robust multiobjective optimization
\newline
\textbf{Mathematics Subject Classification (2020)}\hspace{3mm}41A65 . 49K99 . 65K10 . 90C29 . 90C46

\footnote{\textsf{$^*$Corresponding author}}

\afterpage{} \fancyhead{} \fancyfoot{} \fancyhead[LE, RO]{\bf\thepage} \fancyhead[LO]{\small Approximate solutions for robust multiobjective optimization programming in Asplund spaces} \fancyhead[RE]{\small Saadati and Oveisiha}

\section{Introduction}\label{Sec1-Intro}
\emph{Robust optimization} investigates the cases in which optimization problems often consider uncertain data due to prediction errors, lack of information, fluctuations, or disturbances \cite{2,44,3}. In particular, in such cases these problems rely on conflicting goals due to different multiobjective optimization criteria. Hence, the \emph{robust multiobjective optimization} is highly of interest in optimization theory and substantial in applications.

The first concept of robustness as a kind of sensitivity against perturbations for multiobjective optimization problems was explored by Branke \cite{42} and provided by Deb and Gupta \cite{10}. In addition, various concepts in minimax robustness for multiobjective optimization were introduced by Kuroiwa and Lee \cite{13}, Jeyakumar et al. \cite{5}, Ehrgott et al. \cite{17}, and Ide and K{\"o}bis \cite{18}. Recently, some different concepts of robustness used in multiobjective optimization in the face of data uncertainty have been established in \cite{6,54,8,43,25}.

Approximate efficient solutions of multiobjective optimization problems can be viewed as feasible points whose objective values display a prescribed error $\varepsilon$ in the optimal values of the vector objective. This concept has been widely studied in \cite{52,53,64,65}. Optimality conditions and duality theories of $\varepsilon\mbox{-}\textrm{efficient}$ solutions and $\varepsilon\mbox{-}\textrm{quasi}\mbox{-}\textrm{efficient}$ solutions for convex programming problems under uncertainty have been presented in \cite{60,61}.

The most significant results have been introduced to approximate robust optimization in the finite-dimensional case. So, an infinite-dimensional framework would be proper to study when involving optimality and duality in approximate robust multiobjective optimization. From this, we are motivated to articulate and analyze problems that consider infinite-dimensional frameworks.

Let \mbox{$f: X \to Y$} be a locally Lipschitzian vector-valued function between \emph{Asplund} spaces, and let $\Omega \subset X$ be a nonempty closed set. Suppose that $K \subset Y$ be a pointed (i.e., $K \bigcap \,(-K) = \{0\}$) closed convex cone. We consider the following \emph{multiobjective} optimization problem:
\begin{equation*}\hypertarget{P}{}
\begin{aligned}
  (\mathrm{P}) \qquad \min\nolimits_{K} \,\,\, &f(x) \\
                        \textrm{s.t.} \,\,\, &g_{i}(x) \le 0, \quad i \in 1,2,\dots,n, \nonumber
\end{aligned}
\end{equation*}
where the functions $g_{i}: X \to \mathbb{R}$, $i = 1,2,\dots,n$, define the constraints. Problem (\hyperlink{P}{P}) under data \emph{uncertainty} in the constraints can be captured by the following \emph{uncertain multiobjective} optimization problem:
\begin{equation*}\hypertarget{UP}{}
\begin{aligned}
  (\mathrm{UP}) \qquad \min\nolimits_{K} \,\,\, &f(x) \\
                        \textrm{s.t.} \,\,\, &g_{i}(x, v_{i}) \le 0, \quad i \in 1,2,\dots,n, \nonumber
\end{aligned}
\end{equation*}
where $x \in X$ is the vector of \emph{decision} variable, $v_{i}$ is the vector of \emph{uncertain} parameter and $v_{i} \in \mathcal{V}_{i}$ for some \emph{sequentially compact} topological space $\mathcal{V}_{i}$, $v := (v_{1},v_{2},\dots,v_{n}) \in \mathcal{V} := \prod\limits_{i=1}^{n} \mathcal{V}_{i}$, and $g_{i}: X \times \mathcal{V}_{i} \to \mathbb{R}$, $i = 1,2,\dots,n$, are given functions.

One of the powerful deterministic structures to study problem (\hyperlink{UP}{UP}) is the \emph{robust} optimization, which is known as the problem that the uncertain objective and constraint are satisfied for all possible scenarios within a prescribed uncertainty set. We now associate with them:
\begin{equation*}\hypertarget{RP}{}
\begin{aligned}
  (\mathrm{RP}) \qquad \min\nolimits_{K} \,\,\, &f(x) \\
                        \textrm{s.t.} \,\,\, &g_{i}(x, v_{i}) \le 0, \quad \forall v_{i} \in \mathcal{V}_{i}, \,\, i = 1,2,\dots,n.
\end{aligned}
\end{equation*}
The feasible set $F$ of problem (\hyperlink{RP}{RP}) is defined by
\begin{equation*}
  F := \big\{ x \in \Omega \,\mid\, g_{i}(x, v_{i}) \le 0, \,\, \forall v_{i} \in \mathcal{V}_{i}, \, i = 1,2,\dots,n \big\}.
\end{equation*}

\begin{defn}\hypertarget{Def2-2}{}
Let $\vartheta \in K$, one says a vector $\bar{x} \in F$ is
\begin{itemize}
  \item [(i)] a \emph{robust} $\vartheta\mbox{-}\mathit{efficient}$ \emph{solution} of problem (\hyperlink{UP}{UP}), denoted by $\bar{x} \in \vartheta\mbox{-}\mathcal{S}(RP)$, iff
      \begin{equation*}
        f(x) - f(\bar{x}) + \vartheta \notin - K \setminus \{0\}, \quad \forall x \in  F,
      \end{equation*}
  \item [(ii)] a \emph{weakly robust} $\vartheta\mbox{-}\mathit{efficient}$ \emph{solution} of problem (\hyperlink{UP}{UP}), denoted by $\bar{x} \in \vartheta\mbox{-}\mathcal{S}^{w}(RP)$, iff
      \begin{equation*}
        f(x) - f(\bar{x}) + \vartheta \notin - \textrm{int}\hspace{.4mm}K, \quad \forall x \in  F,
      \end{equation*}
  \item [(iii)] a \emph{robust} $\vartheta\mbox{-}\mathit{quasi}\mbox{-}\mathit{efficient}$ \emph{solution} of problem (\hyperlink{UP}{UP}), denoted by $\bar{x} \in \vartheta\mbox{-}\textrm{quasi}\mbox{-}\mathcal{S}(RP)$, iff
      \begin{equation*}
        f(x) - f(\bar{x}) + \|x - \bar{x}\| \, \vartheta \notin - K \setminus \{0\}, \quad \forall x \in  F,
      \end{equation*}
  \item [(iv)] a \emph{weakly robust} $\vartheta\mbox{-}\mathit{quasi}\mbox{-}\mathit{efficient}$ \emph{solution} of problem (\hyperlink{UP}{UP}), denoted by $\bar{x} \in \vartheta\mbox{-}\textrm{quasi}\mbox{-}\mathcal{S}^{w}(RP)$, iff
      \begin{equation*}
        f(x) - f(\bar{x}) + \|x - \bar{x}\| \, \vartheta \notin - \textrm{int}\hspace{.4mm}K, \quad \forall x \in  F.
      \end{equation*}
\end{itemize}
\end{defn}

The organization of this paper is as follows. In Section \ref{Sec2-Preli}, we recall some preliminary definitions from variational analysis and several auxiliary results. Section \ref{Sec3-NecSuf} provides necessary condition for weakly robust $\vartheta\mbox{-}\textrm{efficient}$ solutions and also necessary/sufficient optimality conditions for (weakly) robust $\vartheta\mbox{-}\textrm{quasi}\mbox{-}\textrm{efficient}$ solutions of problem (\hyperlink{UP}{UP}) in the sense of the limiting subdifferential. In Section \ref{Sec4-Duality}, we formulate duality relations for (weakly) robust $\vartheta\mbox{-}\textrm{quasi}\mbox{-}\textrm{efficient}$ solutions between the corresponding problems.

\section{Preliminaries}\label{Sec2-Preli}
Our notation and terminology are basically standard in the area of variational analysis; see, e.g., \cite{27}. Throughout this paper, all the spaces are Asplund, unless otherwise stated, with the norm $\|\cdot\|$ and the canonical pairing $\langle \cdot\,,\cdot \rangle$ between the space $X$ in question and its \emph{dual} $X^{*}$ equipped with the \emph{weak$^{*}$ topology $w^{*}$}. By $B_{X}(x,r)$, we denote the \emph{closed ball} centered at $x \in X$ with radius $r > 0$, while $B_{X}$ and $B_{X^{*}}$ stand for the \emph{closed unit ball} in $X$ and $X^{*}$, respectively. Given a nonempty set $ \Omega \subset X $, the symbols $\textrm{co}\hspace{.4mm}\Omega$, $\textrm{cl}\hspace{.4mm}\Omega$, and $\textrm{int}\hspace{.4mm}\Omega$ signify the \emph{convex hull}, \emph{topological closure}, and \emph{topological interior} of $\Omega$, respectively, while $\textrm{cl}^{*}\Omega$ stands for the \emph{weak$^{*}$ topological closure} of $\Omega \subset  X^{*}$. The \emph{dual cone} of $\Omega$ is the set
\begin{equation*}
  \Omega^{+} := \big\{ x^{*} \in X^{*} \,\mid\, \langle x^{*}, x \rangle \ge 0, \,\,\, \forall x \in \Omega \big\}.
\end{equation*}
Furthermore, $\mathbb{R}^{n}_{+}$ indicates the nonnegative orthant of $\mathbb{R}^{n}$ for $n \in \mathbb{N} := \{1,2,\dots\}$.

A given set-valued mapping $H : \Omega \subset X \overrightarrow{\to} X^{*}$ is called \emph{weak$^{*}$ closed} at $\bar{x} \in \Omega$ if for any sequence $\{x_{k}\} \subset \Omega$, $x_{k} \to \bar{x}$, and any sequence $\{x^{*}_{k}\} \subset X^{*}$, $x^{*}_{k} \in H(x_{k})$, $x^{*}_{k} \overset{w^{*}} \to x^{*}$, one has $x^{*} \in H(\bar{x})$.

For a set-valued mapping $H : X \overrightarrow{\to} X^{*}$, the \emph{sequential Painlev\'{e}-Kuratowski upper/outer limit} of $H$ as $x \to \bar{x}$ is defined by
\begin{align*}
  \underset{x \to \bar{x}}{\textrm{Lim}\sup} \, H(x) := \Big\{ x^{*} \in X^{*} \,\mid\,\,\, &\exists
  \text{ sequences } x_{k} \to \bar{x} \text{ and } x^{*}_{k} \overset{{\scriptscriptstyle w^{*}}} \to x^{*} \\
  &\text{with } x^{*}_{k} \in H(x_{k}) \text{ for all } k \in \mathbb{N} \Big\}.
\end{align*}

Let $\Omega \subset X$ be \emph{locally closed} around $\bar{x} \in \Omega$, i.e., there is a neighborhood $U$ of $\bar{x}$ for which $\Omega \bigcap \textrm{cl}\hspace{.4mm}U$ is closed. The \emph{Fr\'{e}chet normal cone} $\widehat{N}(\bar{x}; \Omega)$ and the \emph{Mordukhovich normal cone} $N(\bar{x}; \Omega)$ to $\Omega$ at $\bar{x} \in \Omega$ are defined by
\begin{align}
  \label{2-1}
  \widehat{N}(\bar{x}; \Omega) &:= \Big\{x^{*} \in X^{*} \,\mid\, \limsup\limits_{x \overset{\hspace{-1mm}\sOmega} \to \bar{x}} \dfrac{\langle x^{*}, x - \bar{x} \rangle}{\|x - \bar{x}\|} \le 0\Big\}, \\
  \label{2-2}
  N(\bar{x}; \Omega) &:= \underset{x \overset{\hspace{-1mm}\sOmega} \to \bar{x}}{\textrm{Lim}\sup} \, \widehat{N}(x; \Omega),
\end{align}
where $x \overset{\hspace{-1mm}\Omega} \to \bar{x}$ stands for $x \to \bar{x}$ with $x \in \Omega$. If $\bar{x} \notin \Omega$, we put $\widehat{N}(\bar{x}; \Omega) = N(\bar{x}; \Omega) := \emptyset$.

For an extended real-valued function $\phi : X \to \overline{\mathbb{R}}$, the \emph{limiting/Mordukhovich subdifferential} and the \emph{regular/Fr\'{e}chet subdifferential} of $\phi$ at $\bar{x} \in \textrm{dom}\,\phi$ are given, respectively, by
\begin{equation*}
  \partial \phi(\bar{x}) := \big\{ x^{*} \in X^{*} \,\mid\, (x^{*}, -1) \in  N((\bar{x}, \phi(x)); \textrm{epi}\,\phi) \big\}
\end{equation*}
and
\begin{equation*}
  \widehat{\partial} \phi(\bar{x}) := \big\{ x^{*} \in X^{*} \,\mid\, (x^{*}, -1) \in  \widehat{N}((\bar{x}, \phi(x)); \textrm{epi}\,\phi) \big\}.
\end{equation*}
If $|\phi(\bar{x})| = \infty$, then one puts $\partial \phi(\bar{x}) := \widehat{\partial} \phi(\bar{x}) := \emptyset$.

Assign $\langle y^{*}, f \rangle (x) := \langle y^{*}, f(x) \rangle$, $x \in X$, $y^{*} \in Y^{*}$, for a vector-valued map $f : X \to Y$, and denote $\textrm{gph}\,f := \big\{ (x, y) \in X \times Y \,\mid\, y = f(x) \big\}$. Next we recall the required results known as the scalarization formulae of the \emph{coderivatives}.

\begin{lem}\label{Lem2-1}
Let $y^{*} \in Y^{*}$, and let $f : X \to Y$ be Lipschitz around $\bar{x} \in X$. We have
\begin{itemize}
  \item [\emph{(i)}] \emph{(See \cite[Proposition~3.5]{28})} $x^{*} \in \widehat{\partial} \langle y^{*}, f \rangle(\bar{x}) \,\, \Leftrightarrow \,\, (x^{*}, -y^{*}) \in \widehat{N}((\bar{x}, f(\bar{x})); \text{gph}\,\,f)$.
  \item [\emph{(ii)}] \emph{(See \cite[Theorem~1.90]{27})} $x^{*} \in \partial \langle y^{*}, f \rangle(\bar{x}) \,\, \Leftrightarrow \,\, (x^{*}, -y^{*}) \in N((\bar{x}, f(\bar{x})); \text{gph}\,\,f)$.
\end{itemize}
\end{lem}

Another calculus result is the \emph{sum rule} for the limiting subdifferential.

\begin{lem}\label{Lem2-2} \emph{(See \cite[Theorem~3.36]{27})}
Let $ \phi_{i} : X \to \overline{\mathbb{R}}$, $(i \in \{1, 2,\dots, n\}, n \ge 2)$, be lower semicontinuous around $\bar{x}$, and let all but one of these functions be Lipschitz continuous around $\bar{x} \in X$. Then, one has
\begin{equation*}
  \partial ( \phi_{1} + \phi_{2} +\dots+ \phi_{n})(\bar{x}) \subset \partial \phi_{1}(\bar{x}) + \partial \phi_{2}(\bar{x}) + \dots + \partial \phi_{n}(\bar{x}).
\end{equation*}
\end{lem}

The following lemma computes the limiting subdifferential for the \emph{maximum} functions in Asplund spaces. The interested readers are referred to \cite{25,50,51} for more details and proofs. The notation $\partial_{x}$ indicates the limiting subdifferential operation with respect to $x$.

\begin{lem}\label{Lem2-6}
Let $\mathcal{V}$ be a sequentially compact topological space, and let $g : X \times \mathcal{V} \to \mathbb{R}$ be a function such that for each fixed $v \in \mathcal{V}$, $g(\cdot, v)$ is locally Lipschitz on $U \subset X$ and for each fixed $x \in U$, $g(x, \cdot)$ is upper semicontinuous on $\mathcal{V}$. Let $\phi(x) := \max\limits_{v \in \mathcal{V}} g(x, v)$. If the multifunction $(x, v) \in U \times \mathcal{V} \,\, \overrightarrow{\to} \,\, \partial_{x} g(x, v) \subset X^{*}$ is weak$^{*}$ closed at $(\bar{x}, \bar{v})$ for each $\bar{v} \in \mathcal{V}(\bar{x})$, then the set \mbox{$\emph{cl}^{*}\emph{co} \Big(\medcup \Big\{\partial_{x} g(\bar{x}, v) \,\mid\, v\in \mathcal{V}(\bar{x})\Big\}\Big)$} is nonempty and
\begin{equation*}
\partial \phi(\bar{x}) \subset \emph{cl}^{*}\emph{co} \Big(\medcup \Big\{\partial_{x} g(\bar{x}, v) \,\mid\, v \in \mathcal{V}(\bar{x})\Big\}\Big),
\end{equation*}
where $\mathcal{V}(\bar{x}) = \big\{ v \in \mathcal{V} \,\mid\, g(\bar{x}, v) = \phi(\bar{x}) \big\} $.
\end{lem}

In what follows, we also use a formula for the limiting subdifferential of maximum of finitely many functions in Asplund spaces.
\begin{lem}\label{Lem2-7}\emph{(See \cite[Theorem~3.46]{27})} Let $\phi_{i} : X \to \overline{\mathbb{R}}$, $(i \in \{1, 2,\dots, n\}, n \ge 2)$, be Lipschitz continuous around $\bar{x}$. Put $\phi(x) := \max\limits_{i \in \{1, 2,\dots, n\}} \phi_{i}(x)$. Then
\begin{equation*}
  \partial \phi(\bar{x}) \subset \medcup \Big\{\partial \Big(\sum_{i \in I(\bar{x})} \mu_{i} \, \phi_{i} \Big)(\bar{x}) \,\mid\, (\mu_{1}, \mu_{2},\dots, \mu_{n}) \in \Lambda(\bar{x})\Big\},
\end{equation*}
where
\begin{equation*}
  I(\bar{x}) := \big\{ i \in \{1, 2,\dots, n\} \,\mid\, \phi_{i}(\bar{x}) = \phi(\bar{x}) \big\}
\end{equation*}
and
\begin{equation*}
  \Lambda(\bar{x}) := \Big\{(\mu_{1}, \mu_{2},\dots, \mu_{n}) \,\mid\, \mu_{i} \ge 0, \,\, \sum_{i=1}^{n} \mu_{i} = 1, \,\, \mu_{i} \, (\phi_{i}(\bar{x}) - \phi(\bar{x})) = 0 \Big\}.
\end{equation*}
\end{lem}

\begin{assum}\hypertarget{assum}{}(\hspace{-.05mm}See \cite[p.131]{25})
Suppose $\mathcal{V}$ be a sequentially compact topological space, and let $f:X \to Y$ and $g: X \times \mathcal{V} \to \mathbb{R}^{n}$ are functions satisfying the following hypotheses:
\begin{itemize}
  \item[(A1)] For a fixed $\bar{x} \in \Omega$, $g$ is locally Lipschitz in the first argument and uniformly on $\mathcal{V}$ in the second argument, i.e., there exist an open neighborhood $U$ of $\bar{x}$ and a positive constant $\ell$ such that $\|g(z, v) - g(y, v)\| \le \ell \|z - y\|$ for all $z, y \in U$ and $v \in \mathcal{V}$.
  \item[(A2)] For each $i = 1, 2, \dots, n$, the function $v_{i} \in \mathcal{V}_{i} \mapsto g_{i}(x, v_{i}) \in \mathbb{R}$ is upper semicontinuous for each $x \in U$.
  \item[(A3)] For each $i = 1, 2, \dots, n$, we define real-valued functions $\phi_{i}$ and $\phi$ on $X$ via
      \begin{equation*}
         \phi_{i}(x) := \max_{v_{i} \in \mathcal{V}_{i}} g_{i}(x, v_{i}) \,\,\,\, \text{ and } \,\,\,\, \phi(x) := \max_{i \in \{1, 2, \dots, n\}} \phi_{i}(x),
      \end{equation*}
      and we notice that above assumptions imply that $\phi_{i}$ is well defined on $\mathcal{V}_{i}$. In addition, $\phi_{i}$ and $\phi$ follow readily that are locally Lipschitz on $U$, since each $g_{i}(\cdot, v_{i})$ is (see \cite[(H1), p.131]{25} and \cite[p.290]{6}). Note that the feasible set $F$ can be equivalently characterized by:
      \begin{equation*}
        F = \big\{x \in \Omega \,\mid\, \phi_{i}(x) \le 0, \,\, i = 1, 2, \dots, n\big\} = \big\{x \in \Omega \,\mid\, \phi(x) \le 0\big\}.
      \end{equation*}
  \item[(A4)] For each $i = 1, 2, \dots, n$, the multifunction $(x, v_{i}) \in U \times \mathcal{V}_{i} \,\, \overrightarrow{\to} \,\, \partial_{x} g_{i}(x, v_{i}) \subset X^{*}$ is weak$^{*}$ closed at $(\bar{x}, \bar{v}_{i})$ for each $\bar{v}_{i} \in \mathcal{V}_{i}(\bar{x})$, where $\mathcal{V}_{i}(\bar{x}) = \big\{ v_{i} \in \mathcal{V}_{i} \,\mid\, g_{i}(\bar{x}, v_{i}) = \phi_{i}(\bar{x}) \big\}$.
   \item[(A5)] For  a fixed $\bar{x} \in \Omega$, $\vartheta \in K$, and $y^{*} \in K^{+}$, we define a real-valued function $\psi$ on $X$ as follows:
      \begin{equation*}
         \psi(x) := \max \big\{\langle y^{*}, f(x) - f(\bar{x}) + \vartheta \rangle, \phi(x) \big\}.
      \end{equation*}
\end{itemize}
\end{assum}

Inspired by the concept of pseudo-quasi generalized convexity by Fakhar \cite{52}, we introduce a similar
concept of pseudo-quasi convexity type for $f$ and $g$.

\begin{defn}\hypertarget{Def2-4}{}
Let $\vartheta \in K$, we say that
\begin{itemize}
\item [(i)] $(f, g)$ is $\vartheta\mbox{-}\mathit{type}$ \emph{I pseudo convex} on $\Omega$ at $\bar{x} \in \Omega$ if for any $x \in \Omega$, $y^{*} \in K^{+}$, $u^{*} \in {\partial} \langle y^{*}, f \rangle(\bar{x})$, and $v^{*}_{i} \in \partial_{x} g_{i}(\bar{x}, v_{i})$, $v_{i} \in \mathcal{V}_{i}(\bar{x})$, $i = 1, 2, \dots, n$, there exists $w \in -N(\bar{x}; \Omega)^{+}$ such that
    \begin{align*}
    &\langle y^{*}, f \rangle(x) < \langle y^{*}, f \rangle(\bar{x}) - \|x - \bar{x}\| \langle y^{*}, \vartheta \rangle \Longrightarrow \langle u^{*}, w \rangle + \|x - \bar{x}\| \langle y^{*}, \vartheta \rangle < 0,\\
    &g_{i}(x, v_{i}) \le g_{i}(\bar{x}, v_{i}) \Longrightarrow \langle v^{*}_{i}, w \rangle \le 0, \quad i = 1, 2, \dots, n,\\
    &\|w\| \le \|x - \bar{x}\|.
    \end{align*}
\item [(ii)] $(f, g)$ is $\vartheta\mbox{-}\mathit{type}$ \emph{II pseudo convex} on $\Omega$ at $\bar{x} \in \Omega$ if for any $x \in \Omega \setminus \{\bar{x}\}$, $y^{*} \in K^{+} \setminus \{0\}$, $u^{*} \in {\partial} \langle y^{*}, f \rangle(\bar{x})$, and $v^{*}_{i} \in \partial_{x} g_{i}(\bar{x}, v_{i})$, $v_{i} \in \mathcal{V}_{i}(\bar{x})$, $i = 1, 2, \dots, n$, there exists $w \in -N(\bar{x}; \Omega)^{+}$ such that
    \begin{align*}
    &\langle y^{*}, f \rangle(x) \le \langle y^{*}, f \rangle(\bar{x}) - \|x - \bar{x}\| \langle y^{*}, \vartheta \rangle \Longrightarrow \langle u^{*}, w \rangle + \|x - \bar{x}\| \langle y^{*}, \vartheta \rangle < 0,\\
    &g_{i}(x, v_{i}) \le g_{i}(\bar{x}, v_{i}) \Longrightarrow \langle v^{*}_{i}, w \rangle \le 0, \quad i = 1, 2, \dots, n,\\
    &\|w\| \le \|x - \bar{x}\|.
    \end{align*}
\end{itemize}
\end{defn}

\begin{rem}\hypertarget{Exa2-2(Rem)}{}
If in Definition \hyperlink{Def2-4}{2.1},
\begin{itemize}
  \item [(i)] we set $\Omega = X$ and $\vartheta = 0$, then this definition reduces to \cite[Definition~2.2]{54}.
  \item [(ii)] we set $Y = \mathbb{R}^{p}$, then this definition reduces to \cite[Definition~3.8]{52}.
  \item [(ii)] we set $\Omega = X$, $Y = \mathbb{R}^{p}$, and $\vartheta = 0$, then this definition reduces to \cite[Definition~3.2]{8}.
\end{itemize}
\end{rem}

\begin{rem}\hypertarget{Exa2-1(Rem)}{}
\begin{itemize}
  \item [(i)] It follows from Definition \hyperlink{Def2-4}{2.1} that if $(f, g)$ is $\vartheta\mbox{-}\textrm{type}$ II pseudo convex on $\Omega$ at $\bar{x} \in \Omega$, then $(f, g)$ is $\vartheta\mbox{-}\textrm{type}$ I pseudo convex on $\Omega$ at $\bar{x} \in \Omega$, but converse is not true (see Example \ref{Exa2-2}).
  \item [(ii)] It is noted that the generalized (resp., strictly generalized) convexity (see \cite[Definition~3.2]{52}) of $(f, g)$ is reduced to the $\vartheta\mbox{-}\textrm{type}$ I (resp., type II) pseudo convexity of $(f, g)$. Furthermore, as the next example demonstrates, the class of $\vartheta\mbox{-}\textrm{type}$ I pseudo convex functions is properly wider than the class of generalized convex functions, which is properly broader than convex functions (see \cite[Example~3.12]{53}).
\end{itemize}
\end{rem}

\begin{exa}\label{Exa2-2}
Let $X := \mathbb{R}^{2}$, $Y := \mathbb{R}^{3}$, $\Omega := \mathbb{R}^{2}$, $\mathcal{V}_{i} := [-1, -\dfrac{1}{4}]$, $i = 1, 2$, $\mathcal{V} := \prod\limits_{i=1}^{2} \mathcal{V}_{i}$, and let $K := \{(y_{1},y_{2},y_{3}) \in \mathbb{R}^{3} \,\mid\, y_{1} \le 0 \text{ and } y_{i} \ge 0 \text{ for } i=2,3 \}$. Consider $f : X \to Y$ and $g : X \times \mathcal{V} \to \mathbb{R}^{2}$ defined by $f := (f_{1}, f_{2}, f_{3})$ and $g := (g_{1}, g_{2})$, respectively, where
\begin{equation*}
        \left\{\begin{aligned}
              f_{1}(x_{1},x_{2}) &= 5 |x_{1}| - \frac{2}{5} x_{2} + \frac{4}{5}, \\
              f_{2}(x_{1},x_{2}) &= \frac{1}{2} |x_{1}| + 6, \\
              f_{3}(x_{1},x_{2}) &= 4 |x_{1}| + \frac{1}{2} x_{2} + 1
              \end{aligned}
        \right.
        \quad \text{and} \quad
        \left\{\begin{aligned}
              g_{1}(x_{1},x_{2}, v_{1}) &= \frac{1}{4} v_{1}^{2} |x_{1}| + \frac{1}{2} v_{1}^{2} x_{2} - v_{1}^{2} + \frac{1}{4} |v_{1}|, \\
              g_{2}(x_{1},x_{2}, v_{2}) &= \frac{1}{8} x_{1}^{2} + |v_{2}| x_{2} - |v_{2}| + \frac{1}{4},
              \end{aligned}
        \right.
\end{equation*}
$v_{i} \in \mathcal{V}_{i}$, $i = 1,2$. Let $\vartheta := (0,0,\dfrac{3}{2}) \in K$ and consider $\bar{x} := (0, 0) \in \Omega$, Hence $N(\bar{x}; \Omega) = \{(0,0)\}$ and $N(\bar{x}; \Omega)^{+} = \mathbb{R}^{2}$. Obviously, from the definitions,
\begin{equation*}
  \partial f_{1}(\bar{x}) = [-5, 5] \times \{-\dfrac{2}{5}\}, \,\,\,\, \partial f_{2}(\bar{x}) = [-\dfrac{1}{2}, \dfrac{1}{2}] \times \{0\}, \,\, \text{ and } \,\, \partial f_{3}(\bar{x}) = [-4, 4] \times \{\dfrac{1}{2}\}.
\end{equation*}
Moreover
\begin{equation*}
  \partial_{x} g_{1}(\bar{x}, v_{1}) = [-\dfrac{1}{4} v_{1}^{2}, \dfrac{1}{4} v_{1}^{2}] \times \{\dfrac{1}{2} v_{1}^{2}\} \,\, \text{ and } \,\, \partial_{x} g_{2}(\bar{x}, v_{2}) = (0, |v_{2}|)
\end{equation*}
for all $v_{i} \in \mathcal{V}_{i}$, $i = 1,2$.

Suppose that for some $x := (x_{1}, x_{2}) \in \Omega$ and $y^{*} := (y^{*}_{1}, y^{*}_{2}, y^{*}_{3}) \in K^{+}$ the condition $\langle y^{*}, f \rangle(x) < \langle y^{*}, f \rangle(\bar{x}) - \|x - \bar{x}\| \langle y^{*}, \vartheta \rangle$ is satisfied. Thus
\begin{equation*}
  (-\frac{2}{5} y^{*}_{1} + \frac{1}{2} y^{*}_{3}) x_{2} < - (5 y^{*}_{1} + \frac{1}{2} y^{*}_{2} + 4 y^{*}_{3}) |x_{1}| - \|x - \bar{x}\| \langle y^{*}, \vartheta \rangle.
\end{equation*}
Dividing both sides of above inequality by $c := -\dfrac{2}{5} y^{*}_{1} + \dfrac{1}{2} y^{*}_{3} > 0$, we have
\begin{equation}\label{2-18}
  x_{2} < - \frac{1}{c} (5 y^{*}_{1} + \frac{1}{2} y^{*}_{2} + 4 y^{*}_{3}) |x_{1}| - \frac{1}{c} \|x - \bar{x}\| \langle y^{*}, \vartheta \rangle.
\end{equation}
Putting $(w_{1},w_{2}) := w = x \in -N(\bar{x}; \Omega)^{+}$ and employing (\ref{2-18}), for any $u^{*} := y^{*}_{1} u^{*}_{1} + y^{*}_{2} u^{*}_{2} + y^{*}_{3} u^{*}_{3} \in \partial \langle y^{*}, f \rangle(\bar{x})$, where $u^{*}_{i} := (u^{*}_{ix}, u^{*}_{iy}) \in \partial f_{i}(\bar{x})$, $i = 1, 2, 3$, one has $\|w\| \le \|x - \bar{x}\|$ and
\begin{align*}
  \langle u^{*}, w \rangle + \|x - \bar{x}\| \langle y^{*}, \vartheta \rangle & = \begin{pmatrix}
                                           y^{*}_{1} u^{*}_{1x} + y^{*}_{2} u^{*}_{2x} + y^{*}_{3} u^{*}_{3x} \\
                                           -\dfrac{2}{5} y^{*}_{1} + \dfrac{1}{2} y^{*}_{3}
                                         \end{pmatrix}
                                         \times
                                         \begin{pmatrix}
                                           w_{1} \\
                                           w_{2}
                                         \end{pmatrix}
                                         + \|x - \bar{x}\| \langle y^{*}, \vartheta \rangle \\
                                     & < (y^{*}_{1} u^{*}_{1x} + y^{*}_{2} u^{*}_{2x} + y^{*}_{3} u^{*}_{3x}) x_{1} - (5 y^{*}_{1} + \frac{1}{2} y^{*}_{2} + 4 y^{*}_{3}) |x_{1}| \\
                                     & \le 0,
\end{align*}
where the latter inequality is due to $u^{*}_{1x} \in [-5, 5]$, $u^{*}_{2x} \in [-\dfrac{1}{2}, \dfrac{1}{2}]$, and $u^{*}_{3x} \in [-4, 4]$. So $\langle u^{*}, w \rangle + \|x - \bar{x}\| \langle y^{*}, \vartheta \rangle < 0$.

Now, for $x \in \Omega$ and $v_{1} \in \mathcal{V}_{1}(\bar{x})$ the condition $g_{1}(x, v_{1}) \le g_{1}(\bar{x}, v_{1})$ implies $\dfrac{1}{2} v_{1}^{2} x_{2} \le -\dfrac{1}{4} v_{1}^{2} |x_{1}|$, and therefore for any $v^{*}_{1} := (v^{*}_{1x}, v^{*}_{1y}) \in \partial_{x}g_{1}(\bar{x}, v_{1})$, we get
\begin{equation*}
  \langle v^{*}_{1}, w \rangle \le v^{*}_{1x} x_{1} - \frac{1}{4} v_{1}^{2} |x_{1}| \le 0
\end{equation*}
due to $v^{*}_{1x} \in [-\dfrac{1}{4} v_{1}^{2}, \dfrac{1}{4} v_{1}^{2}]$. Similarly for $x \in \Omega$ and $v_{2} \in \mathcal{V}_{2}(\bar{x})$ satisfying $g_{2}(x, v_{2}) \le g_{2}(\bar{x}, v_{2})$, we have $x_{2} \le - \dfrac{1}{8|v_{1}|} x_{1}^{2}$, and thus for any $v^{*}_{2} \in \partial_{x}g_{2}(\bar{x}, v_{2})$, it holds
\begin{equation*}
  \langle v^{*}_{2}, w \rangle \le - \dfrac{1}{8} x_{1}^{2} \le 0.
\end{equation*}
Therefore, $(f, g)$ is $\vartheta\mbox{-}\textrm{type}$ I pseudo convex on $\Omega$ at $\bar{x}$.

Although, there exist $x := (0, 1) \in \Omega \setminus \{\bar{x}\}$ and $y^{*} := (0, 1, 0) \in K^{+} \setminus \{0\}$ such that $\langle y^{*}, f \rangle (x) = 6 = \langle y^{*}, f \rangle (\bar{x}) - \|x - \bar{x}\| \langle y^{*}, \vartheta \rangle$, but for $u_{1}^{*} := (0, -\dfrac{2}{5}) \in \partial f_{1}(\bar{x})$, $u_{2}^{*} := (0,0) \in \partial f_{2}(\bar{x})$, and $u_{3}^{*} := (0, \dfrac{1}{2}) \in \partial f_{3}(\bar{x})$, one has
\begin{equation*}
  u^{*} := y_{1}^{*} u_{1}^{*} + y_{2}^{*} u_{2}^{*} + y_{3}^{*} u_{3}^{*} = (0,0)
\end{equation*}
so $\langle u^{*}, w \rangle + \|x - \bar{x}\| \langle y^{*}, \vartheta \rangle = 0$ for any $w \in -N(\bar{x}; \Omega)^{+}$. This signifies that $(f, g)$ is not $\vartheta\mbox{-}\textrm{type}$ II pseudo convex on $\Omega$ at $\bar{x}$. On the other side, there exist $x := (-1, -5) \in \Omega$ and $y^{*} := (0, 0, 0) \in K^{+}$ such that for any $w \in -N(\bar{x}; \Omega)^{+}$ with $\|w\| \le \|x - \bar{x}\|$ we have
\begin{align*}
  \langle y^{*}, f \rangle (x) &- \langle y^{*}, f \rangle (\bar{x}) = 0, \\
  g_{1}(x, v_{1}) &- g_{1}(\bar{x}, v_{1}) = - \dfrac{9}{4} v_{1}^{2} < 0, \\
  g_{2}(x, v_{2}) &- g_{2}(\bar{x}, v_{2}) = \dfrac{1}{8} - 5 |v_{2}| < 0.
\end{align*}
Hence, $(f, g)$ is not generalized convex on $\Omega$ at $\bar{x}$.
\end{exa}

\begin{exa}\label{Exa2-3}
Suppose $X$, $Y$, $\Omega$, $\mathcal{V}_{i}$, $i = 1,2$, $\mathcal{V} := \prod\limits_{i=1}^{2} \mathcal{V}_{i}$, and $K$ be the same as Example \ref{Exa2-2}. Let $f : X \to Y$ defined by $f := (f_{1}, f_{2}, f_{3})$, where
\begin{equation*}
        \left\{\begin{aligned}
              f_{1}(x_{1}, x_{2}) &= -\frac{4}{5} x_{1}^{2} + 5|x_{1}| - \frac{4}{5} x_{2}^{2} - \frac{2}{5} x_{2} + \frac{4}{5}, \\
              f_{2}(x_{1}, x_{2}) &= \frac{1}{2} |x_{1}| + 6, \\
              f_{3}(x_{1}, x_{2}) &= x_{1}^{2} + 4 |x_{1}| + x_{2}^{2} + \frac{1}{2} x_{2} + 1,
              \end{aligned}
        \right.
\end{equation*}
and let $g : X \times \mathcal{V} \to \mathbb{R}^{2}$ be the same as Example \ref{Exa2-2}. Let $\bar{x} := (0, 0) \in \Omega$ and $\vartheta$ be the same as Example \ref{Exa2-2}. Then
\begin{equation*}
  \partial f_{1}(\bar{x}) = [-5, 5] \times \{-\dfrac{2}{5}\}, \,\,\,\, \partial f_{2}(\bar{x}) = [-\dfrac{1}{2}, \dfrac{1}{2}] \times \{0\}, \,\, \text{ and } \,\, \partial f_{3}(\bar{x}) = [-4, 4] \times \{\dfrac{1}{2}\}.
\end{equation*}

Suppose that for some $x := (x_{1}, x_{2}) \in \Omega \setminus \{\bar{x}\}$ and $y^{*} := (y^{*}_{1}, y^{*}_{2}, y^{*}_{3}) \in K^{+} \setminus \{0\}$ the condition $\langle y^{*}, f \rangle(x) \le \langle y^{*}, f \rangle(\bar{x}) - \|x - \bar{x}\| \langle y^{*}, \vartheta \rangle$ is satisfied. Therefore
\begin{align*}
  (-\frac{2}{5} y^{*}_{1} + \frac{1}{2} y^{*}_{3}) x_{2} \le &- (-\frac{4}{5} y^{*}_{1} + y^{*}_{3}) x_{1}^{2} - (5 y^{*}_{1} + \frac{1}{2} y^{*}_{2} + 4 y^{*}_{3}) |x_{1}| \\
  &- (-\frac{4}{5} y^{*}_{1} + y^{*}_{3}) x_{2}^{2} - \|x - \bar{x}\| \langle y^{*}, \vartheta \rangle.
\end{align*}
Dividing both sides of above inequality by $c := \dfrac{2}{5} y^{*}_{1} + \dfrac{1}{2} y^{*}_{3} > 0$, we have
\begin{equation}\label{2-22}
  x_{2} \le - 2 x_{1}^{2} - \dfrac{1}{c} (5 y^{*}_{1} + \dfrac{1}{2} y^{*}_{2} + 4 y^{*}_{3}) |x_{1}| - 2 x_{2}^{2} - \dfrac{1}{c} \|x - \bar{x}\| \langle y^{*}, \vartheta \rangle.
\end{equation}
Putting $(w_{1},w_{2}) := w = x \in -N(\bar{x}; \Omega)^{+}$ and using (\ref{2-22}), for any $u^{*} := y^{*}_{1} u^{*}_{1} + y^{*}_{2} u^{*}_{2} + y^{*}_{3} u^{*}_{3} \in \partial \langle y^{*}, f \rangle(\bar{x})$, where $u^{*}_{i} := (u^{*}_{ix}, u^{*}_{iy}) \in \partial f_{i}(\bar{x})$, $i = 1, 2, 3$, we get
\begin{align*}
  \langle u^{*}, w \rangle + \|x - \bar{x}\| \langle y^{*}, \vartheta \rangle & \le (y^{*}_{1} u^{*}_{1x} + y^{*}_{2} u^{*}_{2x} + y^{*}_{3} u^{*}_{3x}) x_{1} - 2 c x_{1}^{2} - (5 y^{*}_{1} + \frac{1}{2} y^{*}_{2} + 4 y^{*}_{3}) |x_{1}| - 2 c x_{2}^{2} \\
 & < 0,
\end{align*}
where the latter strict inequality is due to $u^{*}_{1x} \in [-5, 5]$, $u^{*}_{2x} \in [-\dfrac{1}{2}, \dfrac{1}{2}]$, $u^{*}_{3x} \in [-4, 4]$, and $x \ne \bar{x}$. So $\langle u^{*}, w \rangle + \|x - \bar{x}\| \langle y^{*}, \vartheta \rangle < 0$. The complete calculation is similar to that of Example \ref{Exa2-2}. Hence, $(f, g)$ is $\vartheta\mbox{-}\textrm{type}$ II pseudo convex on $\Omega$ at $\bar{x}$.
\end{exa}

In the rest of this section, we present a suitable constraint qualification in the sense of robust, which is required to get a so-called \emph{robust} $\vartheta\mbox{-}\mathit{approximate}$ \emph{Karush-Kuhn-Tucker} (KKT) \emph{condition}.

\begin{defn}\hypertarget{Def2-5}{}(\hspace{-.05mm}See \cite[Definition 4.5]{52})
Let $\bar{x} \in F$. We say that the \emph{constraint qualification} (CQ) \emph{condition} is satisfied at $\bar{x}$ if
\begin{equation*}
  0 \notin \textrm{cl}^{*}\textrm{co}\Big(\medcup \Big\{\partial_{x} g_{i}(\bar{x}, v_{i}) \,\mid\, v_{i} \in \mathcal{V}_{i}(\bar{x})\Big\}\Big) + N(\bar{x}; \Omega),\quad i \in I(\bar{x}),
\end{equation*}
where $I(\bar{x}) := \big\{ i \in \{1, 2,\dots, n\} \,\mid\, \phi_{i}(\bar{x}) = \phi(\bar{x}) \big\}$.
\end{defn}

It is noteworthy here that this condition (CQ) is reduced to the \emph{extended Mangasarian-Fromovitz constraint qualification} (EMFCQ) in the \emph{smooth} setting; see e.g., \cite{27} for more details.

\begin{defn}\hypertarget{Def2-6}{}
Let $\vartheta \in K$ and $\bar{x} \in F$. One says that $\bar{x}$ satisfies the \emph{robust} $\vartheta\mbox{-}\mathit{approximate}$ (KKT) \emph{condition} if there exist $y^{*} \in K^{+} \setminus \{0\}$, $\mu := (\mu_{1}, \mu_{2},\dots, \mu_{n}) \in \mathbb{R}^{n}_{+}$, and $\bar{v}_{i} \in \mathcal{V}_{i}$, $i=1,2,\dots,n$, such that
\begin{equation*}
 \left\{\begin{aligned}
      & 0 \in \partial \langle y^{*}, f \rangle(\bar{x}) + \sum_{i = 1}^{n} \mu_{i} \, \textrm{cl}^{*}\textrm{co} \Big(\medcup \Big\{\partial_{x} g_{i}(\bar{x}, v_{i}) \,\mid\, v_{i} \in \mathcal{V}_{i}(\bar{x})\Big\}\Big) + \langle y^{*}, \vartheta \rangle B_{X^{*}} + N(\bar{x}; \Omega), \\
      & \mu_{i} \, \max\limits_{v_{i} \in \mathcal{V}_{i}} g_{i}(\bar{x}, v_{i}) = \mu_{i} \, g_{i}(\bar{x}, \bar{v}_{i}) = 0, \quad i = 1, 2, \dots, n.
       \end{aligned}
 \right.
\end{equation*}
\end{defn}

Therefore, the robust $\vartheta\mbox{-}\textrm{approximate}$ (KKT) condition defined above is guaranteed by the constraint qualification (CQ).

\section{Robust necessary and sufficient optimality conditions}\label{Sec3-NecSuf}
This section is devoted to study necessary optimality conditions for weakly robust $\vartheta\mbox{-}\textrm{efficient}$ solutions and weakly robust $\vartheta\mbox{-}\textrm{quasi}\mbox{-}\textrm{efficient}$ solutions of problem (\hyperlink{UP}{UP}) by exploiting the nonsmooth version of Fermat's rule, the sum rule for the limiting subdifferential and the scalarization formulae of the coderivatives, and to discuss sufficient optimality conditions for (weakly) robust $\vartheta\mbox{-}\textrm{quasi}\mbox{-}\textrm{efficient}$ solutions by imposing the pseudo convexity assumptions.

The first theorem presents a necessary optimality condition in a \emph{fuzzy form} for weakly robust $\vartheta\mbox{-}\textrm{efficient}$ solutions of problem (\hyperlink{UP}{UP}).

\begin{thm}\label{Thm3-3}
Suppose that $f$ and $g_{i}$, $i = 1, 2, \dots, n$, satisfy \emph{\bf Assumptions} \hyperlink{assum}{\emph{(A1)-(A5)}}. If $\bar{x} \in \vartheta\mbox{-}\mathcal{S}^{w}(RP)$, then there exist $x_{\eta} \in \Omega$, with $\|x_{\eta} - \bar{x}\| \le \eta$, $y^{*} \in K^{+}$, $\lambda := (\lambda_{1},\lambda_{2}) \in \mathbb{R}^{2}_{+}$, with $\lambda_{2} \ne 0$, and $\mu := (\mu_{1},\mu_{2},\dots,\mu_{n}) \in \mathbb{R}^{n}_{+}$, with $ \dfrac{\lambda_{1}}{\lambda_{2}}\|y^{*}\| + \|\mu\| = 1$, and $v_{i_{\eta}} \in \mathcal{V}_{i}$, $i = 1, 2, \dots, n$, such that
\begin{align*}
& 0 \in \dfrac{\lambda_{1}}{\lambda_{2}}\,\partial \langle y^{*}, f \rangle(x_{\eta}) + (1-\lambda_{1}) \sum_{i=1}^{n} \mu_{i} \, \emph{cl}^{*}\emph{co} \Big(\medcup \Big\{\partial_{x} g_{i}(x_{\eta}, v_{i}) \,\mid\, v_{i} \in \mathcal{V}_{i}(x_{\eta})\Big\}\Big) \\
&\hspace{2mm}+ \dfrac{\langle y^{*}, \vartheta \rangle}{\lambda_{2}\,\eta}B_{X^{*}} + N(x_{\eta}; \Omega), \\
& \dfrac{\lambda_{1}}{\lambda_{2}}\Big(\langle y^{*}, f(x_{\eta}) - f(\bar{x}) + \vartheta \rangle - \psi(x_{\eta})\Big) = 0, \\
& (1-\lambda_{1})\mu_{i}\Big(g_{i}(x_{\eta}, v_{i_{\eta}}) - \psi(x_{\eta})\Big) = 0, \quad i = 1, 2, \dots, n.
\end{align*}
\end{thm}

\begin{proof}
If $\bar{x} \in \vartheta\mbox{-}\mathcal{S}^{w}(RP)$, then we have $f(x) - f(\bar{x}) + \vartheta \notin -\textrm{int}\hspace{.4mm}K^{+}$ for all $x \in F$. Using the separation theorem, there exists $y^{*} \in K^{+}$ such that
\begin{equation}\label{3-28}
   \langle y^{*}, f(x) - f(\bar{x}) + \vartheta \rangle \ge 0, \quad \forall x \in  F.
\end{equation}

Let us consider the function $\psi$ and take into account (\ref{3-28}), it can be easily obtained that
\begin{equation}\label{3-23}
0 \le \psi(x), \quad \forall x \in \Omega,
\end{equation}
which implies that $\psi$ is bounded from below on $\Omega$, too.
\\
Furthermore, due to $\bar{x} \in F$, it holds that $\psi(\bar{x}) = \langle y^{*}, \vartheta \rangle$. Thus, from (\ref{3-23}) we get that
\begin{equation*}
\psi(\bar{x}) \le \inf_{x \in \Omega} \psi(x) + \langle y^{*}, \vartheta \rangle.
\end{equation*}

For any $\eta > 0$, using the \emph{Ekeland's variable principle} (see \cite[Theorem~1.1]{62}), we arrive at $x_{\eta} \in \Omega$ such that $\|x_{\eta} - \bar{x}\| \le \eta$ and
\begin{equation*}
\psi(x_{\eta}) \le \psi(x) + \dfrac{\langle y^{*}, \vartheta \rangle}{\eta}\|x_{\eta} - x\|, \quad \forall x \in \Omega.
\end{equation*}
This denotes that $x_{\eta}$ is a minimizer to the optimization problem
\begin{equation*}
\min_{x \in \Omega} \omega(x),
\end{equation*}
where
\begin{equation}\label{3-25}
\omega(x) := \psi(x) + \dfrac{\langle y^{*}, \vartheta \rangle}{\eta}\|x_{\eta} - x\|, \quad x \in \Omega.
\end{equation}
Thus $x_{\eta}$ is a minimizer to the unconstrained optimization problem
\begin{equation}\label{3-26}
\min_{x \in X} \omega(x) + \delta(x;\Omega).
\end{equation}

Applying the \emph{generalized Fermat's rule} (see \cite[Proposition~1.114]{27}), we obtain
\begin{equation}\label{3-29}
  0 \in \partial \big(\omega + \delta(\cdot;\Omega)\big)(x_{\eta}).
\end{equation}
Since the function $\omega$ is Lipschitz continuous around $x_{\eta}$ and the function $\delta(\cdot;\Omega)$ is l.s.c around this point, from the sum rule of Lemma \ref{Lem2-2} applied to (\ref{3-29}) and from the relation $\partial\delta(x_{\eta};\Omega) = N(x_{\eta};\Omega)$ we get
that
\begin{equation}\label{3-30}
0 \in \partial \omega(x_{\eta}) + N(x_{\eta};\Omega).
\end{equation}

Also note that (see \cite[Example~4]{63})
\begin{equation*}
\partial\big(\|x_{\eta} - x\|\big)(x_{\eta}) = B_{X^{*}}.
\end{equation*}
Use the summation rule again to $\omega$ defined in (\ref{3-25}) and using (\ref{3-30}), we arrive at
\begin{equation}\label{3-32}
0 \in \partial \psi(x_{\eta}) + \dfrac{\langle y^{*}, \vartheta \rangle}{\eta}B_{X^{*}} + N(x_{\eta};\Omega).
\end{equation}

Now, applying the formula for the limiting subdifferential of maximum functions in Lemma \ref{Lem2-7}, one has
\begin{align*}
  \partial \psi(x_{\eta}) \subset \medcup \Big\{\partial \big(\alpha_{1} \, \langle y^{*}, f(\cdot) &- f(\bar{x}) + \vartheta \rangle + \alpha_{2} \, \phi(\cdot)\big)(x_{\eta}) \,\mid \, \alpha_{1},\alpha_{2} \ge 0, \,\, \alpha_{1} + \alpha_{2} = 1, \\
  &\alpha_{1} \big(\langle y^{*}, f(x_{\eta}) - f(\bar{x}) + \vartheta \rangle - \psi(x_{\eta})\big) = 0,\,\,
  \alpha_{2} \big(\phi(x_{\eta}) - \psi(x_{\eta}) \big) = 0 \Big\}.
\end{align*}
This together with (\ref{3-32}) and using the sum rule give us $(\bar{\alpha}_{1}, \bar{\alpha}_{2}) \in \mathbb{R}^{2}_{+}$ with $\bar{\alpha}_{1} + \bar{\alpha}_{2} = 1$, such that
\begin{align}
\label{3-41}
&\bar{\alpha}_{1} \big(\langle y^{*}, f(x_{\eta}) - f(\bar{x}) + \vartheta \rangle - \psi(x_{\eta})\big) = 0,\\
\label{3-42}
&\bar{\alpha}_{2} \big(\phi(x_{\eta}) - \psi(x_{\eta}) \big) = 0, \\
\intertext{and}
\label{3-12}
&0 \in \bar{\alpha}_{1} \, \partial \langle y^{*}, f \rangle(x_{\eta}) + \bar{\alpha}_{2} \, \partial \phi(x_{\eta}) + \dfrac{\langle y^{*}, \vartheta \rangle}{\eta}B_{X^{*}} + N(x_{\eta};\Omega).
\end{align}

Invoking again Lemma \ref{Lem2-7}, we have
\begin{equation}\label{3-21}
  \partial \phi(x_{\eta}) \subset \medcup \Big\{\partial \Big(\sum_{i \in I(x_{\eta})} \mu_{i} \, \phi_{i}\Big)(x_{\eta}) \,\mid\, (\mu_{1},\mu_{2},\dots,\mu_{n})\in \Lambda(x_{\eta})\Big\},
\end{equation}
where $I(x_{\eta}) = \big\{ i \in \{1, 2, \dots, n\} \,\mid\, \phi_{i}(x_{\eta}) = \phi(x_{\eta}) \big\}$ and
\begin{equation*}
  \Lambda(x_{\eta}) = \Big\{ (\mu_{1},\mu_{2},\dots,\mu_{n}) \,\mid\, \mu_{i} \ge 0, \,\, \sum_{i=1}^{n} \mu_{i} = 1, \,\, \mu_{i} \, (\phi_{i}(x_{\eta}) - \phi(x_{\eta})) = 0 \Big\}.
\end{equation*}
Using further Lemma \ref{Lem2-6}, we arrive at
\begin{equation}\label{3-22}
  \partial \phi_{i}(x_{\eta}) \subset \textrm{cl}^{*}\textrm{co} \Big(\medcup \Big\{\partial_{x} g_{i}(x_{\eta}, v_{i}) \,\mid\, v_{i} \in \mathcal{V}_{i}(x_{\eta})\Big\}\Big), \quad i = 1, 2, \dots, n,
\end{equation}
where $\mathcal{V}_{i}(x_{\eta}) = \big\{ v_{i} \in \mathcal{V}_{i} \,\mid\, g_{i}(x_{\eta}, v_{i}) = \phi_{i}(x_{\eta}) \big\}$ and the set $ \textrm{cl}^{*}\textrm{co} \Big(\medcup \Big\{\partial_{x} g_{i}(x_{\eta}, v_{i}) \,\mid\, v_{i} \in \mathcal{V}_{i}(x_{\eta})\Big\}\Big)$ is nonempty. The sum rule of the limiting subdifferential and the relations (\ref{3-12})-(\ref{3-22}) results to
\begin{align*}
  0 &\in \bar{\alpha}_{1} \, \partial \langle y^{*}, f \rangle(x_{\eta}) +  \bar{\alpha}_{2} \, \medcup \Big\{\sum_{i\in I(x_{\eta})} \mu_{i} \, \textrm{cl}^{*}\textrm{co} \Big(\medcup \Big\{\partial_{x} g_{i}(x_{\eta}, v_{i}) \,\mid\, v_{i} \in \mathcal{V}_{i}(x_{\eta})\Big\}\Big) \,\mid\, \\
  &(\mu_{1},\mu_{2},\dots,\mu_{n}) \in \Lambda(x_{\eta})\Big\} + \dfrac{\langle y^{*}, \vartheta \rangle}{\eta}B_{X^{*}} + N(x_{\eta}; \Omega).
\end{align*}
So, there exist $\bar{\mu} := (\bar{\mu}_{1},\bar{\mu}_{2},\dots,\bar{\mu}_{n}) \in \Lambda(x_{\eta})$, with $\mathlarger{\sum\limits}_{i=1}^{n} \, \bar{\mu}_{i} = 1$ and $\bar{\mu}_{i} = 0$ for all $i \in \{1, 2, \dots, n\} \setminus I(x_{\eta})$, such that
\begin{equation*}
 0 \in \bar{\alpha}_{1} \, \partial \langle y^{*}, f \rangle(x_{\eta}) + \bar{\alpha}_{2} \, \sum_{i=1}^{n} \bar{\mu}_{i} \, \textrm{cl}^{*}\textrm{co} \Big(\medcup \Big\{\partial_{x} g_{i}(x_{\eta}, v_{i}) \,\mid\, v_{i} \in \mathcal{V}_{i}(x_{\eta})\Big\}\Big) + \dfrac{\langle y^{*}, \vartheta \rangle}{\eta}B_{X^{*}} + N(x_{\eta};\Omega).
\end{equation*}
Putting $\lambda_{1} := \bar{\alpha}_{1}$, $\lambda_{2} := \bar{\alpha}_{1}\,\|y^{*}\| + \|\bar{\mu}\|$, and $\mu := \dfrac{1}{\lambda_{2}} \, \bar{\mu}$ and dividing the above inclusion by $\lambda_{2}$, we have $y^{*} \in K^{+}$, $\lambda := (\lambda_{1},\lambda_{2}) \in \mathbb{R}^{2}_{+}$, $\lambda_{2} \ne 0$, and $\mu \in \mathbb{R}^{n}_{+}$, with $\dfrac{\lambda_{1}}{\lambda_{2}}\| y^{*}\| + \|\mu\| = 1 $, satisfying the first relation in the theorem.

On the other hand, due to \mbox{the upper semicontinuity of the} function $v_{i} \in \mathcal{V}_{i} \longmapsto g_{i}(x_{\eta}, v_{i})$ for each $i = 1, 2, \dots, n$ and also, the sequentially compactness of $\mathcal{V}_{i}$, we can select $v_{i_{\eta}} \in \mathcal{V}_{i}$ such that $g_{i}(x_{\eta}, v_{i_{\eta}}) = \max\limits_{v_{i} \in \mathcal{V}_{i}} g_{i}(x_{\eta}, v_{i}) = \phi_{i}(x_{\eta})$. By considering $\phi_{i}(x_{\eta}) = \phi(x_{\eta})$ for all $i \in I(x_{\eta})$ and using relations (\ref{3-41}) and (\ref{3-42}), the proof of the theorem is achieved.
\end{proof}

\begin{rem}\hypertarget{Exa3-1(Rem)}{}
Theorem \ref{Thm3-3} develops \cite[Theorem~3.4]{53}, where the underlying optimization problem has a finite dimension framework.
\end{rem}

Similarly, we establish a necessary optimality condition in the sense of the limiting subdifferential for weakly robust $\vartheta\mbox{-}\textrm{quasi}\mbox{-}\textrm{efficient}$ solutions of problem (\hyperlink{UP}{UP}). To prove this theorem, it is required to state a fuzzy necessary optimality condition in terms of the Fr\'{e}chet subdifferential for weakly robust $\vartheta\mbox{-}\textrm{quasi}\mbox{-}\textrm{efficient}$ solutions of problem (\hyperlink{UP}{UP}) as follows.

\begin{thm}\label{Thm3-1}\emph{(See \cite[Theorem~3.2]{34})}
Let $\bar{x} \in \vartheta\mbox{-}\textrm{quasi}\mbox{-}\mathcal{S}^{w}(RP)$. Then for each $k \in \mathbb{N}$, there exist $x^{1k} \in B_{X}(\bar{x}, \frac{1}{k})$, \mbox{$x^{2k} \in B_{X}(\bar{x}, \frac{1}{k})$,} $x^{3k} \in \Omega \medcap B_{X}(\bar{x}, \frac{1}{k})$, $y^{*}_{k} \in K^{+}$ with $\| y^{*}_{k} \| = 1$, and $\alpha_{k} \in \mathbb{R}_{+}$ such that
\begin{align*}
   & 0 \in \widehat{\partial} \langle y^{*}_{k}, f \rangle (x^{1k}) + \alpha_{k} \, \widehat{\partial} \phi (x^{2k}) +
   \widehat{N}(x^{3k}; \Omega) + \Big(\langle y^{*}_{k}, \vartheta \rangle + \dfrac{1}{k}\Big) B_{X^{*}}, \nonumber \\
   & | \alpha_{k}\,\phi(x^{2k}) | \le \dfrac{1}{k}.
\end{align*}
\end{thm}

\begin{thm}\label{Thm3-2}
Suppose that $g_{i}$, $i = 1, 2, \dots, n$, satisfy \emph{\bf Assumptions} \hyperlink{assum}{\emph{(A1)-(A4)}}. If $\bar{x} \in \vartheta\mbox{-}\textrm{quasi}\mbox{-}\mathcal{S}^{w}(RP)$, then there exist $y^{*} \in K^{+}$, $\mu := (\mu_{1},\mu_{2},\dots,\mu_{n}) \in \mathbb{R}^{n}_{+}$, with $ \|y^{*}\| + \|\mu\| = 1$, and $\bar{v}_{i} \in \mathcal{V}_{i}$, $i = 1, 2, \dots, n$, such that
\begin{equation}\label{3-13}
 \left\{\begin{aligned}
      & 0 \in \partial \langle y^{*}, f \rangle(\bar{x}) + \sum_{i = 1}^{n} \mu_{i} \, \emph{cl}^{*}\emph{co} \Big(\medcup \Big\{\partial_{x} g_{i}(\bar{x}, v_{i}) \,\mid\, v_{i} \in \mathcal{V}_{i}(\bar{x})\Big\}\Big) + \langle y^{*}, \vartheta \rangle B_{X^{*}} + N(\bar{x}; \Omega), \\
      & \mu_{i} \, \max\limits_{v_{i} \in \mathcal{V}_{i}} g_{i}(\bar{x}, v_{i}) = \mu_{i} \, g_{i}(\bar{x}, \bar{v}_{i}) = 0, \quad i = 1, 2, \dots, n.
       \end{aligned}
 \right.
\end{equation}
Furthermore, if the \emph{(CQ)} is satisfied at $\bar{x}$, then \emph{(\ref{3-13})} holds with $y^{*} \ne 0$.
\end{thm}

\begin{proof}
Let $\bar{x} \in \vartheta\mbox{-}\textrm{quasi}\mbox{-}\mathcal{S}^{w}(RP)$. By using Theorem \ref{Thm3-1}, we obtain sequences $x^{1k} \to \bar{x}$, $x^{2k} \to \bar{x}$, $x^{3k} \to \bar{x}$, $y^{*}_{k} \in K^{+}$ with $\|y^{*}_{k}\| = 1$, $\alpha_{k} \in \mathbb{R}_{+}$, $x^{*}_{1k} \in \widehat{\partial} \langle y^{*}_{k}, f \rangle(x^{1k})$, $x^{*}_{2k} \in \alpha_{k} \, \widehat{\partial} \phi(x^{2k})$, and $x^{*}_{3k} \in \widehat{N}(x^{3k}; \Omega)$ satisfying
\begin{align}
  \label{3-14}
  &0 \in x^{*}_{1k} + x^{*}_{2k} + x^{*}_{3k} + \Big(\langle y^{*}_{k}, \vartheta \rangle + \dfrac{1}{k}\Big) B_{X^{*}}, \\
  &\alpha_{k} \, \phi(x^{2k}) \to 0 \text{ as } k \to \infty. \nonumber
\end{align}
Now we can consider two possibilities for the sequence $\{\alpha_{k}\}$:

\textbf{Case 1:} Suppose that $\{\alpha_{k}\}$ is bounded, therefore without loss of generality we can assume that \mbox{$\alpha_{k} \to \alpha \in \mathbb{R}_{+}$} as $k \to \infty$. In addition, since the sequence $\{y^{*}_{k}\} \subset K^{+}$ is bounded, by applying the weak* sequential compactness of bounded sets in duals to Asplund spaces, there is no loss of generality in assuming that $y^{*}_{k} \overset{w^{*}} \rightarrow \bar{y}^{*} \in K^{+}$ with $\|\bar{y}^{*}\| = 1$ as $k \to \infty$. Let $\ell_{1} > 0$ be a Lipschitz modulus of $f$ around $\bar{x}$. It is obvious that $\|x^{*}_{1k}\| \le \ell_{1} \, \|y^{*}_{k}\| \le \ell_{1}$ for all $k \in \mathbb{N}$ (see, \cite[Proposition~1.85]{27}). As above, by taking a subsequence, if necessary, that $x^{*}_{1k} \overset{w^{*}} \rightarrow x^{*}_{1} \in X^{*}$ as $k \to \infty$. Due to the boundedness of $\{\alpha_{k}\}$ and the Lipschitz continuity of $\phi$ around $\bar{x}$, the sequence $\{x^{*}_{2k}\}$ is also bounded. In this regard, we can have $x^{*}_{2} \in X^{*}$ such that $x^{*}_{2k} \overset{w^{*}} \rightarrow x^{*}_{2} \in X^{*}$ as $k \to \infty$. Using the part (i) of Lemma \ref{Lem2-1} to the inclusion $x^{*}_{1k} \in \widehat{\partial} \langle y^{*}_{k}, f \rangle(x^{1k})$ gives us,
\begin{equation*}
  (x^{*}_{1k}, -y^{*}_{k}) \in \widehat{N}((x^{1k}, f(x^{1k})); \text{gph}\,\,f), \quad k \in \mathbb{N}.
\end{equation*}
Passing the limit as $k \to \infty$ and applying the definitions of normal cones (\ref{2-1}) and (\ref{2-2}), we obtain $(x^{*}_{1}, -\bar{y}^{*}) \in N((\bar{x}, f(\bar{x})); \text{gph}\,\,f)$, which equals to
\begin{equation}\label{3-16}
  x^{*}_{1} \in \partial \langle \bar{y}^{*}, f \rangle(\bar{x}),
\end{equation}
due to the part (ii) of Lemma \ref{Lem2-1}. Similarly, we get
\begin{equation}\label{3-15}
x^{*}_{2} \in \alpha \, \partial \phi(\bar{x}).
\end{equation}

From (\ref{3-14}), there exists $b^{*}_{k} \in B_{X^{*}}$ such that
\begin{equation}\label{3-17}
-x^{*}_{1k} - x^{*}_{2k} - \Big(\langle y^{*}_{k}, \vartheta \rangle + \dfrac{1}{k}\Big) b^{*}_{k} = x^{*}_{3k} \in \widehat{N}(x^{3k}; \Omega), \quad k \in \mathbb{N}.
\end{equation}
Supposing $b^{*}_{k} \to b^{*} \in B_{X^{*}}$ as $k \to \infty$ and passing (\ref{3-17}) to the limit as $k \to \infty$, as well as considering (\ref{2-1}) and (\ref{2-2}), we arrive at
\begin{equation*}
-x^{*}_{1} - x^{*}_{2} - \langle \bar{y}^{*}, \vartheta \rangle b^{*} \in N(\bar{x}; \Omega).
\end{equation*}
Combining the latter with (\ref{3-16}) and (\ref{3-15}) gives us
\begin{equation*}
  0 \in \partial \langle \bar{y}^{*}, f \rangle(\bar{x}) + \alpha \, \partial \phi(\bar{x}) + \langle \bar{y}^{*}, \vartheta \rangle B_{X^{*}} + N(\bar{x}; \Omega).
\end{equation*}

Now, similar to the proof of Theorem \ref{Thm3-3}, there exists $\bar{\mu} := (\bar{\mu}_{1},\bar{\mu}_{2},\dots,\bar{\mu}_{n}) \in \Lambda(\bar{x})$, with $\mathlarger{\sum\limits}_{i=1}^{n} \, \bar{\mu}_{i} = 1$ and $\bar{\mu}_{i} = 0$ for all $i \in \{1, 2, \dots, n\} \setminus I(\bar{x})$, such that
\begin{equation*}
 0 \in \partial \langle \bar{y}^{*}, f \rangle(\bar{x}) + \alpha \sum_{i=1}^{n} \bar{\mu}_{i} \, \textrm{cl}^{*}\textrm{co} \Big(\medcup \Big\{\partial_{x} g_{i}(\bar{x}, v_{i}) \,\mid\, v_{i} \in \mathcal{V}_{i}(\bar{x})\Big\}\Big) + \langle \bar{y}^{*}, \vartheta \rangle B_{X^{*}} + N(\bar{x}; \Omega).
\end{equation*}
Dividing the above inclusion by $\beta := \|\bar{y}^{*}\| + \alpha \, \|\bar{\mu}\|$, and then setting $y^{*} := \dfrac{\bar{y}^{*}}{\beta}$ and $\mu := \dfrac{\alpha}{\beta} \,\bar{\mu}$, we have some $y^{*} \in K^{+}$ and $\mu := (\mu_{1}, \mu_{2},\dots, \mu_{n}) \in \mathbb{R}^{n}_{+}$, with $\|y^{*}\| + \|\mu\| = 1 $, such that
\begin{equation} \label{3-24}
 0 \in \partial \langle y^{*}, f \rangle(\bar{x}) + \sum_{i=1}^{n} \mu_{i} \, \textrm{cl}^{*}\textrm{co} \Big(\medcup \Big\{\partial_{x} g_{i}(\bar{x}, v_{i}) \,\mid\, v_{i} \in \mathcal{V}_{i}(\bar{x})\Big\}\Big) + \langle y^{*}, \vartheta \rangle B_{X^{*}} + N(\bar{x}; \Omega).
\end{equation}

In addition to that, we can obtain $\bar{v}_{i} \in \mathcal{V}_{i}$ such that $g_{i}(\bar{x}, \bar{v}_{i}) = \max\limits_{v_{i} \in \mathcal{V}_{i}} g_{i}(\bar{x}, v_{i}) = \phi_{i}(\bar{x})$ for each $i = 1, 2, \dots, n$. Furthermore, $\alpha \, \phi(\bar{x}) = 0$  since $\alpha_{k} \, \phi(x^{2k}) \to 0$ as $k \to \infty$. By considering that $\phi_{i}(\bar{x}) = \phi(\bar{x})$ for all $i \in I(\bar{x})$, we arrive at
\begin{equation*}
\mu_{i} \, g_{i}(\bar{x}, \bar{v}_{i}) = \dfrac{\alpha}{\beta} \, \bar{\mu}_{i} \, \phi_{i}(\bar{x}) =
\dfrac{\bar{\mu}_{i}}{\beta} \, [\alpha \, \phi(\bar{x})] = 0,
\end{equation*}
i.e., $\mu_{i} \, g_{i}(\bar{x}, \bar{v}_{i}) = \mu_{i} \, \max\limits_{v_{i} \in \mathcal{V}_{i}} g_{i}(\bar{x}, v_{i}) = 0$ for all $i \in \{1, 2, \dots, n\}$. This together with (\ref{3-24}) yields (\ref{3-13}).

\textbf{Case 2:} Next we assume that $\{\alpha_{k}\}$ is unbounded. If $\ell_{2} > 0$ be a Lipschitz constant of $\phi$ around $\bar{x}$, then we have
\begin{equation*}
\|x^{*}_{2k}\| \le \ell_{2} \, \alpha_{k} \text{ \,\,for all } k \in \mathbb{N}.
\end{equation*}

Applying the latter inequality and considering the weak* sequential compactness of bounded sets in duals to Asplund spaces, we may suppose that $\dfrac{x^{*}_{2k}}{\alpha_{k}} \overset{w^{*}} \rightarrow x^{*}_{2} \in X^{*}$ as $k \to \infty$. From (\ref{3-14}), there exists $b^{*}_{k} \in B_{X^{*}}$ such that
\begin{equation}\label{3-18}
-\dfrac{x^{*}_{1k}}{\alpha_{k}} - \dfrac{x^{*}_{2k}}{\alpha_{k}} - \dfrac{\Big(\langle y^{*}_{k}, \vartheta \rangle + \dfrac{1}{k}\Big)b^{*}_{k}}{\alpha_{k}} = \dfrac{x^{*}_{3k}}{\alpha_{k}} \in \widehat{N}(x^{3k}; \Omega), \quad k \in \mathbb{N}.
\end{equation}
Passing (\ref{3-18}) to the limit as $k \to \infty$ and taking (\ref{2-1}) and (\ref{2-2}) into account, we obtain
\begin{equation}\label{3-19}
- x^{*}_{2} \in N(\bar{x}; \Omega).
\end{equation}

Similar to the Case 1, we get from the inclusion $x^{*}_{2k} \in \alpha_{k} \, \widehat{\partial} \phi(x^{2k})$ that
\begin{equation*}
(x^{*}_{2k}, - \alpha_{k}) \in \widehat{N}((x^{2k}, \phi(x^{2k})); \text{gph}\,\,\phi)
\end{equation*}
for each $k \in \mathbb{N}$. So
\begin{equation*}
\Big(\dfrac{x^{*}_{2k}}{\alpha_{k}}, - 1\Big) \in \widehat{N}((x^{2k}, \phi(x^{2k})); \text{gph}\,\,\phi), \quad k \in \mathbb{N}.
\end{equation*}
Assuming $k \to \infty$ and considering (\ref{2-2}) again, we have $(x^{*}_{2}, -1) \in N((\bar{x}, \phi(\bar{x})); \text{gph}\,\,\phi)$, which is equivalent to
\begin{equation*}
x^{*}_{2} \in \partial \phi(\bar{x}).
\end{equation*}
The latter inclusion with (\ref{3-19}) indicate that
\begin{equation*}
0 \in \partial \phi(\bar{x}) + N(\bar{x}; \Omega).
\end{equation*}

Proceed as in the proof of Theorem \ref{Thm3-3}, we have $\mu := (\mu_{1},\mu_{2},\dots,\mu_{n}) \in \mathbb{R}^{n}_{+} \setminus \{0\}$, with $\|\mu\| = 1 $, satisfying
\begin{equation*}
 0 \in \sum_{i=1}^{n} \mu_{i} \, \textrm{cl}^{*}\textrm{co} \Big(\medcup \Big\{\partial_{x} g_{i}(\bar{x}, v_{i}) \,\mid\, v_{i} \in \mathcal{V}_{i}(\bar{x})\Big\} \Big) + N(\bar{x}; \Omega).
\end{equation*}
Moreover, due to the unboundedness of $\{\alpha_{k}\}$ and $\alpha_{k} \, \phi(x^{2k}) \to 0$ as $k \to \infty$, we can choose $\bar{v}_{i} \in \mathcal{V}_{i}$ such that $\mu_{i} \, g_{i}(\bar{x}, \bar{v}_{i}) = \mu_{i} \, \phi_{i}(\bar{x}) =  \mu_{i} \, \phi(\bar{x}) = 0$ for each $i = 1, 2, \dots, n$. So, (\ref{3-13}) holds by taking $y^{*} := 0 \in K^{+}$.

Finally, assuming that $\bar{x}$ satisfies the (CQ) in the Case 1, directly from (\ref{3-13}) we can arrive at $y^{*} \ne 0$, which supports the last statement of the theorem and completes the proof.
\end{proof}

\begin{rem}\hypertarget{Exa3-1(Rem)}{}
Theorem \ref{Thm3-2} reduces to \cite[Theorem~3.2]{54} with $\vartheta = 0$, \cite[Theorem~4.3]{52} with $Y = \mathbb{R}^{p}$, and \cite[Theorem~3.3]{25} and \cite[Theorem~3.7]{53} in the case of finite-dimensional optimization. Note further that our approach here, which involves the fuzzy necessary optimality condition in the sense of the Fr\'{e}chet subdifferential and the inclusion formula for the limiting subdifferential of maximum functions in the setting of Asplund spaces, is totally different from those ones presented in the aforementioned papers.
\end{rem}

\begin{figure}
    \centering
    \includegraphics[width=0.52\textwidth]{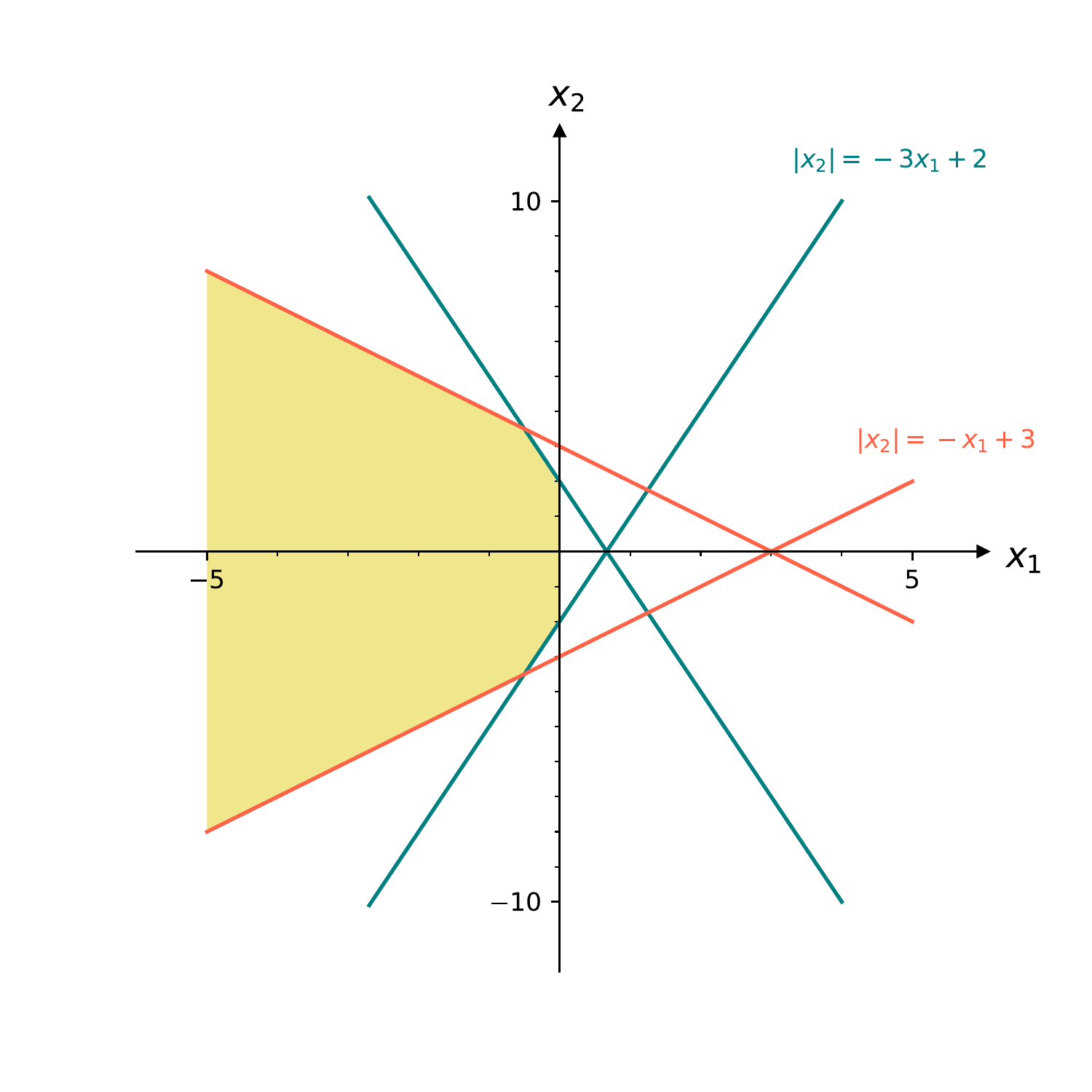}
    \vspace{-1cm}
    \caption{Robust feasible set of problem (UP) in Example \ref{Exa3-2}}
    \label{fig1}
\end{figure}

We then return to an example to illustrate Theorem \ref{Thm3-2} for an uncertain multiobjective optimization problem.
\begin{exa}\label{Exa3-2}
Suppose that $X := \mathbb{R}^{2}$, $Y := \mathbb{R}^{3}$, $\Omega := \mathbb{R}^{2}$, $\mathcal{V}_{i} := [-1, 1]$, $i = 1,2$, $\mathcal{V} := \prod\limits_{i=1}^{2} \mathcal{V}_{i}$, and $K := \mathbb{R}^{3}_{+}$. consider the following uncertain optimization problem:
\begin{equation*}\hypertarget{UPexa3.2}{}
 (\mathrm{UP}) \qquad \min\nolimits_{K} \,\,\, \big\{ f(x) \,\mid\, g_{i}(x, v_{i}) \le 0, \,\, \forall v_{i} \in \mathcal{V}_{i}, \, i = 1,2 \big\},
\end{equation*}
where $f : X \to Y$, $f := (f_{1}, f_{2}, f_{3})$ are given by
\begin{equation*}
        \left\{\begin{aligned}
              f_{1} (x_{1},x_{2}) &= -2 x_{1} + |x_{2}|,\\
              f_{2}(x_{1},x_{2}) &= \dfrac{1}{|x_{1}| + 1} - 3 x_{2} + 2,\\
              f_{3} (x_{1},x_{2}) &= \dfrac{1}{\sqrt{|x_{1}| + 1}} - |x_{2} - 1| - 1,
              \end{aligned}
        \right.
\end{equation*}
and $g : X \times \mathcal{V} \to \mathbb{R}^{2}$, $g := (g_{1}, g_{2})$ are defined by
\begin{equation*}
        \left\{\begin{aligned}
              g_{1}(x_{1},x_{2},v_{1}) &= v_{1}^{2} |x_{2}| + \max \big\{ x_{1}, 2 x_{1} \big\} - 3 |v_{1}|,\\
              g_{2}(x_{1},x_{2},v_{2}) &= -3 |x_{1}| + v_{2} x_{2} - 2,
              \end{aligned}
        \right.
\end{equation*}
where $v_{i} \in \mathcal{V}_{i}$, $i = 1,2$. It is obvious that
\begin{equation*}
  \big\{ v_{1}^{2} |x_{2}| + \max \big\{ x_{1}, 2 x_{1} \big\} - 3 |v_{1}| \le 0 \,\,\, \forall v_{1} \in \mathcal{V}_{1} \big\} = \{(x_{1},x_{2}) \in X \,\mid\, x_{1} \le 0 \text{ and } |x_{2}| \le - x_{1} + 3\},
\end{equation*}
and, due to $x_{1} \le 0$, it can be verified that
\begin{equation*}
  \big\{ -3 |x_{1}| + v_{2} x_{2} - 2 \le 0 \,\,\, \forall v_{2} \in \mathcal{V}_{2} \big\} = \big\{ (x_{1},x_{2}) \in X \,\mid\, x_{1} \le 0 \text{ and } |x_{2}| \le -3 x_{1} + 2 \big\}.
\end{equation*}
Therefore, the robust feasible set is
\begin{align*}
  F = &\big\{(x_{1}, x_{2}) \in X \,\mid\, -\dfrac{1}{2} \le x_{1} \le 0 \text{ and } |x_{2}| \le -3 x_{1} + 2 \big\} \medcup \\
  &\big\{(x_{1}, x_{2}) \in X \,\mid\, x_{1} \le -\dfrac{1}{2} \text{ and } |x_{2}| \le - x_{1} + 3 \big\},
\end{align*}
which is represented in Figure \ref{fig1}.

Let $\vartheta := (\vartheta_{1}, \vartheta_{2}, \vartheta_{3}) = (0,1,0) \in K$ and consider $\bar{x} := (0, 0) \in F$, hence $N(\bar{x}, \Omega) = \{(0,0)\}$. Suppose $x := (x_{1}, x_{2}) \in F$ and take $x_{1} \le 0$, we get $f_{1}(x) - f_{1}(\bar{x}) + \| x - \bar{x} \| \vartheta_{1} \ge 0$. Therefore
\begin{equation*}
f(x) - f(\bar{x}) + \| x - \bar{x} \| \vartheta \notin - \text{int}\hspace{.4mm}K
\end{equation*}
for all $x \in F$, i.e., $\bar{x}$ is a weakly robust $\vartheta\mbox{-}\textrm{quasi}\mbox{-}\textrm{efficient}$ solution of problem (\hyperlink{UPexa3.2}{UP}). Note further that
\begin{align*}
  \phi_{1}(\bar{x}) &= \max\limits_{v_{1} \in \mathcal{V}_{1}} g_{1}(\bar{x}, v_{1}) = \max\limits_{v_{1} \in \mathcal{V}_{1}} (- 3 |v_{1}|) = 0, \\
  \phi_{2}(\bar{x}) &= \max\limits_{v_{2} \in \mathcal{V}_{2}} g_{2}(\bar{x}, v_{2}) = \max\limits_{v_{2} \in \mathcal{V}_{2}} (v_{2} - 2) = -1.
\end{align*}
So $\phi(\bar{x}) = \max \big\{ \phi_{1}(\bar{x}), \phi_{2}(\bar{x}) \big\} = 0$, $\mathcal{V}_{1}(\bar{x}) = \{0\}$, and $\mathcal{V}_{2}(\bar{x}) = \{1\}$. After calculations, we get
\begin{equation*}
\partial f_{1}(\bar{x}) = \{-2\} \times [-1, 1], \,\,\,\, \partial f_{2}(\bar{x}) = [-1, 1] \times \{-3\}, \,\,\,\, \partial f_{3}(\bar{x}) = \big[-\dfrac{1}{2}, \dfrac{1}{2}\big] \times \{-1, 1\},
\end{equation*}
and also
\begin{equation}\label{3-31}
        \left\{\begin{aligned}
              \textrm{cl}^{*}\textrm{co}\Big(\partial_{x}g_{1}(\bar{x}, v_{1}=0)\Big) &= [1, 2] \times \{0\},\\
              \textrm{cl}^{*}\textrm{co}\Big(\partial_{x}g_{2}(\bar{x}, v_{2}=1)\Big) &= [-3, 3] \times \{1\}.
              \end{aligned}
        \right.
\end{equation}
On the other hand, since $I(\bar{x}) = \big\{ i \in \{1, 2\} \,\mid\, \phi_{i}(\bar{x}) = \phi(\bar{x}) \big\} = \{1\}$, it easily follows from (\ref{3-31}) that the (CQ) is satisfied at $\bar{x}$.

Finally, there exist $y^{*} = (\dfrac{\sqrt{2}}{4}, 0, \dfrac{\sqrt{2}}{4}) \in K^{+}$ and $\mu = (\dfrac{1}{2}, 0) \in \mathbb{R}^{2}_{+}$, with $\|y^{*}\| + \|\mu\| = 1$, $b^{*} = (0, 0) \in B_{X^{*}}$, and $a^{*} = (0,0) \in N(\bar{x}, \Omega)$ such that
\begin{align*}
  0 &=
  \begin{pmatrix}
    \dfrac{\sqrt{2}}{4} & 0 & \dfrac{\sqrt{2}}{4}
  \end{pmatrix}
  \begin{pmatrix}
    -2 & 0 & 0 \\
    1 & -3 & -1
  \end{pmatrix}
  +
  \begin{pmatrix}
    \dfrac{1}{2} & 0
  \end{pmatrix}
  \begin{pmatrix}
    \sqrt{2} & 0 \\
    0 & 1
  \end{pmatrix}
  \\
  &+
  \begin{pmatrix}
    \dfrac{\sqrt{2}}{4} & 0  & \dfrac{\sqrt{2}}{4}
  \end{pmatrix}
  \begin{pmatrix}
    0 & 1  & 0
  \end{pmatrix}
  \begin{pmatrix}
    0 & 0
  \end{pmatrix}
  +
  \begin{pmatrix}
    0 & 0
  \end{pmatrix},
\end{align*}
and $\mu_{i} \, \max\limits_{v_{i} \in \mathcal{V}_{i}} g_{i}(\bar{x}, v_{i}) = 0$ for $i = 1, 2$.
\end{exa}

The next theorem establishes a robust $\vartheta\mbox{-}\textrm{approximate}$ (KKT) sufficient optimality condition for (weakly) robust $\vartheta\mbox{-}\textrm{quasi}\mbox{-}\textrm{efficient}$ solutions of problem (\hyperlink{UP}{UP}).

\begin{thm}\label{Thm3-4}
Assume that $\bar{x} \in F$ satisfies the robust $\vartheta\mbox{-}\mathit{approximate}$ \emph{(KKT)} condition.
\begin{itemize}
  \item [\emph{(i)}] If $(f, g)$ is $\vartheta\mbox{-}\mathit{type}$ I pseudo convex on $\Omega$ at $\bar{x}$, then $\bar{x} \in \vartheta\mbox{-}\textrm{quasi}\mbox{-}\mathcal{S}^{w}(\hyperlink{RP}{\emph{RP}})$.
  \item [\emph{(ii)}] If $(f, g)$ is $\vartheta\mbox{-}\mathit{type}$ II pseudo convex on $\Omega$ at $\bar{x}$, then $\bar{x} \in \vartheta\mbox{-}\textrm{quasi}\mbox{-}\mathcal{S}(\hyperlink{RP}{\emph{RP}})$.
\end{itemize}
\end{thm}

\begin{proof}
Let $\bar{x} \in F$ satisfy the robust $\vartheta\mbox{-}\text{approximate}$ (KKT) condition. Therefore, there exist $y^{*} \in K^{+} \setminus \{0\}$, $u^{*} \in \partial \langle y^{*}, f \rangle(\bar{x})$, $\mu_{i} \ge 0$, $v_{i}^{*} \in \textrm{cl}^{*}\textrm{co} \Big(\medcup \Big\{\partial_{x} g_{i}(\bar{x}, v_{i}) \,\mid\, v_{i} \in \mathcal{V}_{i}(\bar{x})\Big\}\Big)$, $i = 1, 2, \dots, n$, and $b^{*} \in B_{X^{*}}$ such that
\begin{align}\label{3-34}
  & -\Big(u^{*} + \sum_{i=1}^{n} \mu_{i} \, v^{*}_{i} + \langle y^{*}, \vartheta \rangle b^{*} \Big) \in N(\bar{x}; \Omega), \\
  \label{3-35}
  & \mu_{i} \, \max_{v_{i} \in \mathcal{V}_{i}} g_{i}(\bar{x}, v_{i}) = 0, \quad i = 1, 2, \dots, n.
\end{align}

Firstly, we justify (i). Argue by contradiction that $\bar{x} \notin \vartheta\mbox{-}\textrm{quasi}\mbox{-}\mathcal{S}^{w}(\hyperlink{RP}{\textrm{RP}})$. Hence, there is $\hat{x} \in F$ such that $f(\hat{x}) - f(\bar{x}) + \|\hat{x} - \bar{x}\| \vartheta \in - \textrm{int}\hspace{.4mm}K$. The latter gives us $\langle y^{*}, f(\hat{x}) - f(\bar{x}) + \|\hat{x} - \bar{x}\| \vartheta \rangle < 0$ (see \cite[Lemma~3.21]{37}). Since $(f, g)$ is $\vartheta\mbox{-}\text{type}$ I pseudo convex on $\Omega$ at $\bar{x}$, we deduce from this inequality that there exists $w \in -N(\bar{x}; \Omega)^{+}$ such that
\begin{equation}\label{3-36}
 \left\{\begin{aligned}
      &\langle u^{*}, w \rangle +  \langle y^{*}, \vartheta \rangle \|\hat{x} - \bar{x}\| < 0, \\
      &\|w\| \le \|\hat{x} - \bar{x}\|.
       \end{aligned}
 \right.
\end{equation}

On the other side, it follows from (\ref{3-34}) for $w$ above that
\begin{equation}\label{3-37}
  \langle u^{*}, w \rangle + \sum_{i=1}^{n} \mu_{i} \, \langle v^{*}_{i}, w \rangle + \langle y^{*}, \vartheta \rangle \langle b^{*}, w \rangle \ge 0.
\end{equation}
The relations (\ref{3-36}) and (\ref{3-37}) entail that
\begin{equation*}
\sum_{i=1}^{n} \mu_{i} \, \langle v^{*}_{i}, w \rangle > 0.
\end{equation*}

To proceed, we assume that there is $i_{0} \in \{1, 2, \dots, n\}$ such that $\mu_{i_{0}} \, \langle v^{*}_{i_{0}}, w \rangle > 0$. Taking into account that $v^{*}_{i_{0}} \in \textrm{cl}^{*}\textrm{co} \Big(\medcup\Big\{\partial_{x} g_{i_{0}}(\bar{x}, v_{i_{0}}) \,\mid\, v_{i_{0}} \in \mathcal{V}_{i_{0}}(\bar{x})\Big\}\Big)$, we get sequence $\{v^{*}_{i_{0}k}\} \subset \textrm{co} \Big(\medcup \Big\{\partial_{x} g_{i_{0}}(\bar{x}, v_{i_{0}}) \,\mid\, v_{i_{0}} \in \mathcal{V}_{i_{0}}(\bar{x})\Big\}\Big)$ such that $v^{*}_{i_{0}k} \overset{{\scriptscriptstyle w^{*}}} \to v^{*}_{i_{0}}$. Hence, due to $\mu_{i_{0}} > 0$, there is $k_{0} \in \mathbb{N}$ such that
\begin{equation}\label{3-38}
  \langle v^{*}_{i_{0}k_{0}}, w \rangle > 0.
\end{equation}
In addition, since $v^{*}_{i_{0}k_{0}} \in \textrm{co} \Big(\medcup \Big\{\partial_{x} g_{i_{0}} (\bar{x},v_{i_{0}}) \,\mid\, v_{i_{0}} \in \mathcal{V}_{i_{0}}(\bar{x})\Big\}\Big)$, there exist $v^{*}_{p} \in \medcup \Big\{\partial_{x} g_{i_{0}}(\bar{x}, v_{i_{0}}) \,\mid\, v_{i_{0}} \in \mathcal{V}_{i_{0}}(\bar{x})\Big\}$ and $\mu_{p} \ge 0$ with $\mathlarger{\sum\limits}_{p=1}^{s} \, \mu_{p} = 1$, $p = 1,2,\dots,s$, $s \in \mathbb{N}$, such that $v^{*}_{i_{0}k_{0}} = \mathlarger{\sum\limits}_{p=1}^{s} \, \mu_{p} \, v^{*}_{p}$. Combining the latter together (\ref{3-38}), we arrive at $\mathlarger{\sum\limits}_{p=1}^{s} \, \mu_{p} \, \langle v^{*}_{p}, w \rangle > 0$. Thus, we can take $p_{0} \in \{1,2,\dots,s\}$ such that
\begin{equation}\label{3-39}
  \langle v^{*}_{p_{0}}, w \rangle > 0,
\end{equation}
and choose $\bar{v}_{i_{0}} \in \mathcal{V}_{i_{0}}(\bar{x})$ satisfying $v^{*}_{p_{0}} \in \partial_{x} g_{i_{0}}(\bar{x}, \bar{v}_{i_{0}})$ due to $v^{*}_{p_{0}} \in \medcup \Big\{\partial_{x} g_{i_{0}}(\bar{x}, v_{i_{0}}) \,\mid\, v_{i_{0}} \in \mathcal{V}_{i_{0}}(\bar{x})\Big\}$. Invoking now definition of type I pseudo convexity of $(f, g)$ on $\Omega$ at $\bar{x}$, we get from (\ref{3-39}) that
\begin{equation}\label{3-40}
  g_{i_{0}}(\hat{x}, \bar{v}_{i_{0}}) > g_{i_{0}}(\bar{x}, \bar{v}_{i_{0}}).
\end{equation}
Note that $\bar{v}_{i_{0}} \in \mathcal{V}_{i_{0}}(\bar{x})$, thus we have $g_{i_{0}}(\bar{x}, \bar{v}_{i_{0}}) = \max\limits_{v_{i_{0}} \in \mathcal{V}_{i_{0}}} g_{i_{0}}(\bar{x}, v_{i_{0}})$ which together with (\ref{3-35}) yields $\mu_{i_{0}} \, g_{i_{0}}(\bar{x}, \bar{v}_{i_{0}}) = 0$. This implies by (\ref{3-40}) that $\mu_{i_{0}} \, g_{i_{0}}(\hat{x}, \bar{v}_{i_{0}}) > 0$, and hence $g_{i_{0}}(\hat{x}, \bar{v}_{i_{0}}) > 0$, which contradicts with the fact that $\hat{x} \in F$ and completes the proof of (i).

Assertion (ii) is proved similarly to the part (i). If $\bar{x} \notin \vartheta\mbox{-}\textrm{quasi}\mbox{-}\mathcal{S}(RP)$, then there exists $\hat{x} \in F$ such that $f(\hat{x}) - f(\bar{x}) + \|\hat{x} - \bar{x}\| \vartheta  \in - K \setminus \{0\}$. Therefore $\hat{x} \ne \bar{x}$ and $\langle y^{*}, f(\hat{x}) - f(\bar{x}) + \|\hat{x} - \bar{x}\| \vartheta \rangle \le 0$. Now by using the definition of type II pseudo convexity of $(f, g)$ on $\Omega$ at $\bar{x}$, we arrive at the result.
\end{proof}

We immediately get the following robust $\vartheta\mbox{-}\textrm{approximate}$ (KKT) sufficient optimality condition from Remark \hyperlink{Exa2-1(Rem)}{2.2}(i) and Theorem \ref{Thm3-4}.
\begin{cor}
Let $\bar{x} \in F$ satisfy the robust $\vartheta\mbox{-}\mathit{approximate}$ \emph{(KKT)} condition and $(f, g)$ is $\vartheta\mbox{-}\mathit{type}$ I pseudo convex on $\Omega$ at $\bar{x}$, then $\bar{x} \in \vartheta\mbox{-}\textrm{quasi}\mbox{-}\mathcal{S}(\hyperlink{RP}{\emph{RP}})$.
\end{cor}

\begin{rem}\hypertarget{Exa3-3(Rem)}{}
Theorem \ref{Thm3-4} reduces to \cite[Theorem~3.4]{54} with $\vartheta = 0$ and \cite[Theorem~4.7]{52} with $Y = \mathbb{R}^{p}$,  and improves \cite[Theorem~3.11]{25}, \cite[Theorem~3.2]{5}, and \cite[Theorem~3.13]{53} under pseudo convexity assumptions.
\end{rem}

Let us present an example to show the viability of our new concept of pseudo convexity for an uncertain multiobjective optimization problem.

\begin{figure}
    \centering
    \includegraphics[width=0.52\textwidth]{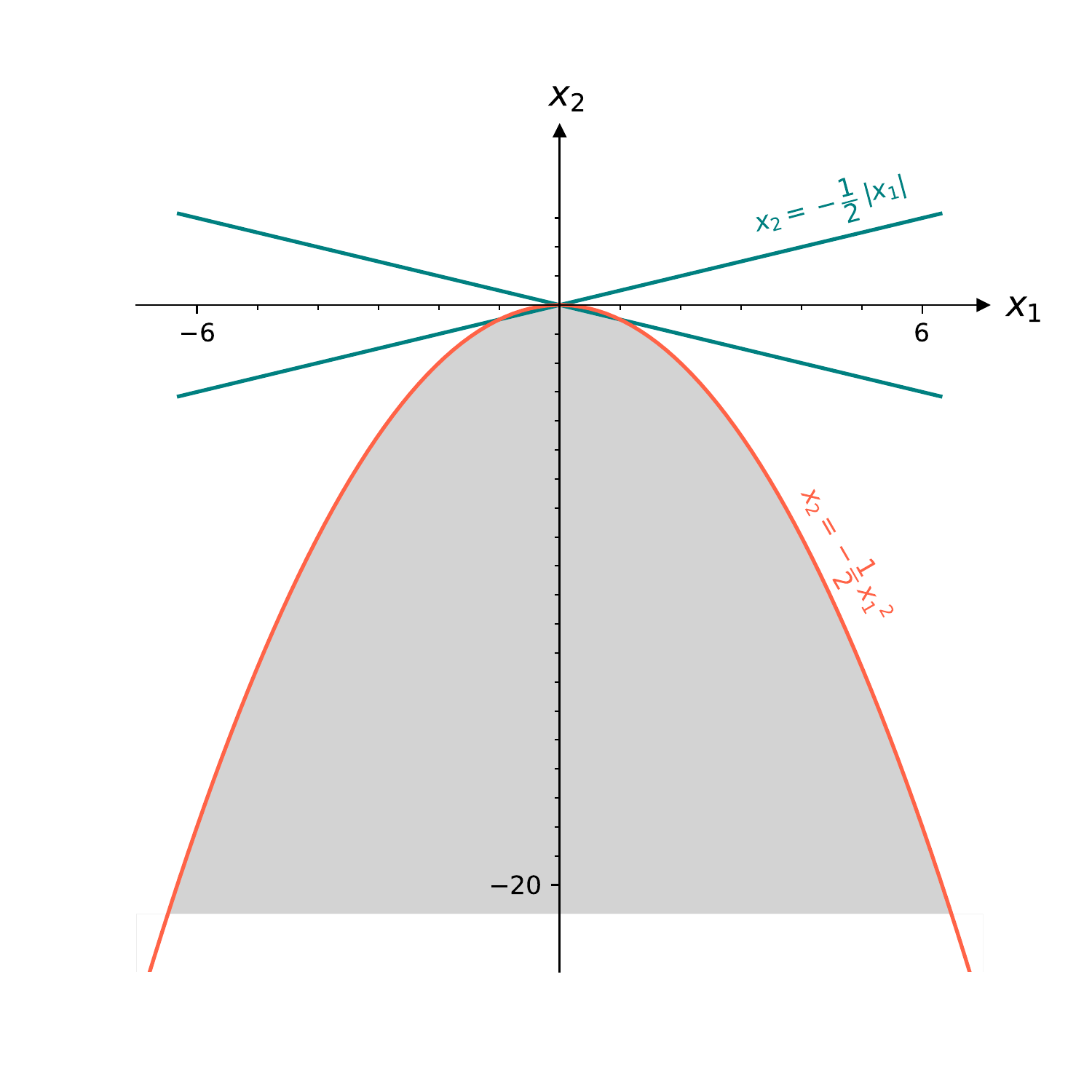}
    \vspace{-1cm}
    \caption{Robust feasible set of problem (UP) in Example \ref{Exa3-5}}
    \label{fig2}
\end{figure}

\begin{exa}\label{Exa3-5}
Let $X$, $Y$, $\Omega$, $\mathcal{V}_{i}$, $i = 1,2$, $\mathcal{V} := \prod\limits_{i=1}^{2} \mathcal{V}_{i}$, $K$, $f$, and $g$ be the same as Example \ref{Exa2-2}. Take the following uncertain optimization problem:
\begin{equation*}\hypertarget{UPexa3.5}{}
 (\mathrm{UP}) \qquad \min\nolimits_{K} \,\,\, \big\{ f(x) \,\mid\, g_{i}(x, v_{i}) \le 0, \,\, \forall v_{i} \in \mathcal{V}_{i}, \, i = 1,2 \big\}.
\end{equation*}
From these constraints of inequality and equality for all $v_{i} \in \mathcal{V}_{i}$, $i = 1,2$, we get
\begin{align*}
  \big\{ \dfrac{1}{4} v_{1}^{2} |x_{1}| + \dfrac{1}{2} v_{1}^{2} x_{2} - v_{1}^{2} + \dfrac{1}{4} |v_{1}| \le 0 \,\,\, \forall v_{1} \in \mathcal{V}_{1} \big\} &= \big\{x \in \mathbb{R}^{2} \,\mid\, x_{2} \le - \dfrac{1}{2} |x_{1}| \big\}, \\
  \big\{ \dfrac{1}{8} x_{1}^{2} + |v_{2}| x_{2} - |v_{2}| + \dfrac{1}{4} \,\,\, \forall v_{2} \in \mathcal{V}_{2} \big\} &= \big\{ x \in \mathbb{R}^{2} \,\mid\, x_{2} \le - \dfrac{1}{2} x_{1}^{2} \big\}.
\end{align*}
So, obviously we can verify that
\begin{equation*}
F= \big\{ x \in \mathbb{R}^{2} \,\mid\, |x_{1}| \le 1 \text{ and } x_{2} \le - \dfrac{1}{2} |x_{1}| \big\} \medcup \big\{ x \in \mathbb{R}^{2} \,\mid\, |x_{1}| > 1 \text{ and } x_{2} \le - \dfrac{1}{2} x_{1}^{2} \big\},
\end{equation*}
as depicted in Figure \ref{fig2}. Let $\bar{x} := (0, 0) \in F$ and $\vartheta := (\vartheta_{1}, \vartheta_{2}, \vartheta_{3}) \in K$ be the same as Example \ref{Exa2-2}. Note that
\begin{align*}
  \phi_{1}(\bar{x}) &= \max\limits_{v_{1} \in \mathcal{V}_{1}} g_{1}(\bar{x}, v_{1}) = \max\limits_{v_{1} \in \mathcal{V}_{1}} (-v_{1}^{2} + \dfrac{1}{4} |v_{1}|) = 0, \\
  \phi_{2}(\bar{x}) &= \max\limits_{v_{2} \in \mathcal{V}_{2}} g_{2}(\bar{x}, v_{2}) = \max\limits_{v_{2} \in \mathcal{V}_{2}} (- |v_{2}| + \dfrac{1}{4}) = 0.
\end{align*}
Then $\phi(\bar{x}) = \max \big\{ \phi_{1}(\bar{x}), \phi_{2}(\bar{x}) \big\} = 0$, $\mathcal{V}_{1}(\bar{x}) = \{-\dfrac{1}{4}\}$, and $\mathcal{V}_{2}(\bar{x}) = \{-\dfrac{1}{4}\}$. It follows from Example \ref{Exa2-2} that $\partial f_{1}(\bar{x}) = [-5, 5] \times \{-\dfrac{2}{5}\}$, $\partial f_{2}(\bar{x}) = [-\dfrac{1}{2}, \dfrac{1}{2}] \times \{0\}$, $\partial f_{3}(\bar{x}) = [-4, 4] \times \{\dfrac{1}{2}\}$, and
\begin{equation*}
\textrm{cl}\hspace{.4mm}\textrm{co} \Big(\partial_{x}g_{1}(\bar{x}, v_{1}=-\dfrac{1}{4})\Big) = [-\dfrac{1}{64}, \dfrac{1}{64}] \times \{\dfrac{1}{32}\}, \quad \textrm{cl}\hspace{.4mm}\textrm{co} \Big(\partial_{x}g_{2}(\bar{x}, v_{2}=-\dfrac{1}{4})\Big) = (0, \dfrac{1}{4}),
\end{equation*}
and further that $(f, g)$ is $\vartheta\mbox{-}\textrm{type}$ I pseudo convex on $\Omega$ at $\bar{x}$. Note that there exist $y^{*} = (-\dfrac{5}{8}, 0, \dfrac{1}{2}) \in K^{+} \setminus \{0\}$, $\mu = (0, 1) \in \mathbb{R}^{2}_{+}$, and $b^{*} = (0, -1) \in B_{X^{*}}$ satisfying
\begin{equation*}
  0 =
  \begin{pmatrix}
    -\dfrac{5}{8} & 0  & \dfrac{1}{2}
  \end{pmatrix}
  \begin{pmatrix}
    \dfrac{8}{5} & 0 & 2 \\
    -\dfrac{2}{5} & 0 & \dfrac{1}{2}
  \end{pmatrix}
  +
  \begin{pmatrix}
    0 & 1
  \end{pmatrix}
  \begin{pmatrix}
    0 & 0 \\
    \dfrac{1}{32} & \dfrac{1}{4}
  \end{pmatrix}
  +
  \begin{pmatrix}
    -\dfrac{5}{8} & 0  & \dfrac{1}{2}
  \end{pmatrix}
  \begin{pmatrix}
    0 & 0  & \dfrac{3}{2}

  \end{pmatrix}
  \begin{pmatrix}
    0 & -1
  \end{pmatrix}
\end{equation*}
and $\mu_{i} \, \max\limits_{v_{i} \in \mathcal{V}_{i}} g_{i}(\bar{x}, v_{i}) = 0$ for $i = 1, 2$. Therefore, the robust $\vartheta\mbox{-}\textrm{approximate}$ (KKT) condition is satisfied at $\bar{x}$. We have
\begin{equation*}
f_{2}(x) - f_{2}(\bar{x}) + \|x - \bar{x}\| \vartheta_{2} = \dfrac{1}{2} |x_{1}| + \|x - \bar{x}\| \vartheta_{1} \ge 0
\end{equation*}
so $f(x) - f(\bar{x}) + \|x - \bar{x}\| \vartheta \notin - \text{int}\hspace{.4mm}K$ for all $x \in F$, i.e., $\bar{x}$ is a weakly robust $\vartheta\mbox{-}\textrm{quasi}\mbox{-}\textrm{efficient}$ solution of problem (\hyperlink{UPexa3.5}{UP}).

On the other hand, if $f$ is the same as Example \ref{Exa2-3}, then $(f, g)$ is $\vartheta\mbox{-}\textrm{type}$ II pseudo convex on $\Omega$ at $\bar{x}$. We get $f_{2}(x) - f_{2}(\bar{x}) + \|x - \bar{x}\| \vartheta_{2} > 0$ for all $x \in F \setminus \{\bar{x}\}$, thus $f(x) - f(\bar{x}) + \|x - \bar{x}\| \vartheta \notin -K \setminus \{0\}$ for all $x \in F$, i.e., $\bar{x}$  is a robust $\vartheta\mbox{-}\textrm{quasi}\mbox{-}\textrm{efficient}$ solution of problem (\hyperlink{UPexa3.5}{UP}).
\end{exa}

\section{Robust duality}\label{Sec4-Duality}
In this section, we formulate the $\vartheta\mbox{-}\mathit{Mond}\mbox{-}\mathit{Weir}\mbox{-}\mathit{type}$ \emph{dual robust problem} (\hyperlink{RD}{$\textrm{RD}_{MW}$}) for (\hyperlink{RP}{RP}), and explore the weak, strong, and converse duality relations between the corresponding problems under pseudo convexity assumptions.

Given $\vartheta \in K$, in connection with the problem (\hyperlink{RP}{RP}), we introduce a \emph{dual robust multiobjective optimization} problem in the sense of Mond-Weir as follows:
\begin{equation*}\hypertarget{RD}{}
  (\textrm{RD}_{MW}) \qquad \max\nolimits_{K} \,\,\,\big\{ \bar{f}(z, y^{*}, \mu) := f(z) \,\mid\, (z, y^{*}, \mu) \in F_{MW} \big\},
\end{equation*}
where $F_{MW}$ is the feasible set defined by
\begin{align*}
  F_{MW} := \bigg\{(z, y^{*}, \mu) \in\,\, &\Omega \times K^{+} \setminus \{0\} \times \mathbb{R}^{n}_{+} \,\mid\, 0 \in \partial \langle y^{*}, f \rangle(z) + \sum_{i=1}^{n} \mu_{i} \, v^{*}_{i} + \langle y^{*}, \vartheta \rangle B_{X^{*}} + N(z; \Omega), \\
  & v_{i}^{*} \in \textrm{cl}^{*}\textrm{co} \Big(\medcup \Big\{\partial_{x} g_{i}(z, v_{i}) \,\mid\, v_{i} \in \mathcal{V}_{i}(z) \Big\}\Big), \,\, \mu_{i} \, g_{i}(z, v_{i}) \ge 0, \,\, i = 1, 2, \dots, n \bigg\}.
\end{align*}

\begin{defn}\hypertarget{Def4-1}{}
Let $\vartheta \in K$, one says a vector $(\bar{z}, \bar{y}^{*}, \bar{\mu}) \in F_{MW}$ is
\begin{itemize}
  \item [(i)] a \emph{robust} $\vartheta\mbox{-}\mathit{quasi}\mbox{-}\mathit{efficient}$ \emph{solution} of problem (\hyperlink{RD}{$\textrm{RD}_{MW}$}), denoted by $(\bar{z}, \bar{y}^{*}, \bar{\mu}) \in \vartheta\mbox{-}\textrm{quasi}\mbox{-}\mathcal{S}(\hyperlink{RD}{\textrm{RD}_{MW}})$, iff
      \begin{equation*}
        \bar{f}(\bar{z}, \bar{y}^{*}, \bar{\mu}) \npreceq \bar{f}(z, y^{*}, \mu) - \|z - \bar{z}\|\vartheta, \quad \forall (z, y^{*}, \mu) \in F_{MW},
      \end{equation*}
  \item [(ii)] a \emph{weakly robust} $\vartheta\mbox{-}\mathit{quasi}\mbox{-}\mathit{efficient}$ \emph{solution} of problem (\hyperlink{RD}{$\textrm{RD}_{MW}$}), denoted by $(\bar{z}, \bar{y}^{*}, \bar{\mu}) \in \vartheta\mbox{-}\textrm{quasi}\mbox{-}\mathcal{S}^{w}(\hyperlink{RD}{\textrm{RD}_{MW}})$, iff
      \begin{equation*}
        \bar{f}(\bar{z}, \bar{y}^{*}, \bar{\mu}) \nprec \bar{f}(z, y^{*}, \mu) - \|z - \bar{z}\|\vartheta, \quad \forall (z, y^{*}, \mu) \in F_{MW}.
      \end{equation*}
  \end{itemize}
\end{defn}

In what follows, we use the following notations for convenience:
\begin{align*}
u \prec v \Leftrightarrow u-v \in -\textrm{int}\hspace{.4mm}K, \quad & u \nprec v \text{ is the negation of } u \prec v, \\
u \preceq v \Leftrightarrow u-v \in -K\setminus\{0\}, \quad & u \npreceq v \text{ is the negation of } u \preceq v.
\end{align*}

Weak duality relations between the primal problem (\hyperlink{RP}{RP}) and the dual problem (\hyperlink{RD}{$\textrm{RD}_{MW}$}) is declared in the following theorem.
\begin{thm}\label{Thm4-1}\textbf{\textsc{(Weak Duality)}}
Let $x \in F$, and let $(z, y^{*}, \mu) \in F_{MW}$.
\begin{itemize}
  \item [\emph{(i)}] If $(f, g)$ is $\vartheta\mbox{-}\mathit{type}$ I pseudo convex on $\Omega$ at $z$, then $f(x) \nprec \bar{f}(z, y^{*}, \mu) - \|x - z\| \vartheta$.
  \item [\emph{(ii)}] If $(f, g)$ is $\vartheta\mbox{-}\mathit{type}$ II pseudo convex on $\Omega$ at $z$, then $f(x) \npreceq \bar{f}(z, y^{*}, \mu) - \|x - z\| \vartheta$.
  \end{itemize}
\end{thm}

\begin{proof}
By $(z, y^{*}, \mu) \in F_{MW}$, there exist $u^{*} \in \partial \langle y^{*}, f \rangle(z)$, $\mu_{i} \ge 0$, $v_{i}^{*} \in \textrm{cl}^{*}\textrm{co} \Big(\medcup \Big\{\partial_{x} g_{i}(z, v_{i}) \,\mid\, v_{i} \in \mathcal{V}_{i}(z)\Big\}\Big)$, $i = 1, 2, \dots, n$, and $b^{*} \in B_{X^{*}}$ such that
\begin{align}\label{4-1}
  & -\Big(u^{*} + \sum_{i=1}^{n} \mu_{i} \, v^{*}_{i} + \langle y^{*}, \vartheta \rangle b^{*} \Big) \in N(z; \Omega),\\
  &  \mu_{i} \, g_{i}(z, v_{i}) \ge 0, \quad i = 1, 2, \dots, n. \nonumber
\end{align}

To prove (i), suppose that $f(x) \prec \bar{f}(z, y^{*}, \mu) - \|x - z\| \vartheta$. Hence $\langle y^{*}, f(x) - \bar{f}(z, y^{*}, \mu) + \|x - z\| \vartheta \rangle < 0$ due to $y^{*} \ne 0$. This is nothing else but $\langle y^{*}, f(x) - f(z) + \|x - z\| \vartheta \rangle < 0$. Since $(f, g)$ is $\vartheta\mbox{-}\textrm{type}$ I pseudo convex on $\Omega$ at $z$, we infer from the last inequality that there exists $w \in -N(z; \Omega)^{+}$ such that
\begin{align*}
&\langle u^{*}, w \rangle +  \langle y^{*}, \vartheta \rangle \|x - z\| < 0, \\
&\|w\| \le \|x - z\|.
\end{align*}

Besides, it follows from (\ref{4-1}) for $w$ above that
\begin{equation*}
  \langle u^{*}, w \rangle + \sum_{i=1}^{n} \mu_{i} \, \langle v^{*}_{i}, w \rangle + \langle y^{*}, \vartheta \rangle \langle b^{*}, w \rangle \ge 0.
\end{equation*}
Combining the latter relations, we get that
\begin{equation*}
  \sum_{i=1}^{n} \mu_{i} \, \langle v^{*}_{i}, w \rangle > 0.
\end{equation*}

Now suppose that there is $i_{0} \in \{1, 2, \dots, n\}$ such that $\mu_{i_{0}} \, \langle v^{*}_{i_{0}}, w \rangle > 0$. Proceeding similarly to the proof of Theorem \ref{Thm3-4}(i) and replacing $\hat{x} - \bar{x}$ with $x - z$ give us $g_{i_{0}}(x, \bar{v}_{i_{0}}) > 0$, which contradicts with $x \in F$.

Next to justify (ii), we proceed similarly to the part (i) by employing $\vartheta\mbox{-}\textrm{type}$ II pseudo convexity of $(f, g)$ on $\Omega$ at $z$, if $f(x) \preceq \bar{f}(z, y^{*}, \mu) - \|x - z\| \vartheta$, then $x \ne z$ and we infer that there exists $w \in -N(z; \Omega)^{+}$ such that $\langle u^{*}, w \rangle +  \langle y^{*}, \vartheta \rangle \|x - z\| < 0$ and $\|w\| \le \|x - z\|$.
\end{proof}

We now establish a strong duality theorem which holds between (\hyperlink{RP}{RP}) and (\hyperlink{RD}{$\textrm{RD}_{MW}$}).

\begin{thm}\label{Thm4-2}\textbf{\textsc{(Strong Duality)}}
Let $\bar{x} \in \mathcal{S}^{w}(\hyperlink{RP}{\emph{RP}})$ be such that the \emph{(CQ)} is satisfied at this point. Then, there exists $(\bar{y}^{*}, \bar{\mu}) \in K^{+} \setminus \{0\} \times \mathbb{R}^{n}_{+}$ such that $(\bar{x}, \bar{y}^{*}, \bar{\mu}) \in F_{MW}$. Furthermore,
\begin{itemize}
  \item [\emph{(i)}] If $(f, g)$ is $\vartheta\mbox{-}\mathit{type}$ I pseudo convex on $\Omega$ at $z$ for all $z \in \Omega$, then $(\bar{x}, \bar{y}^{*}, \bar{\mu}) \in \mathcal{S}^{w}(\hyperlink{RD}{\emph{RD}_{MW}})$.
  \item [\emph{(ii)}] If $(f, g)$ is $\vartheta\mbox{-}\mathit{type}$ II pseudo convex on $\Omega$ at $z$ for all $z \in \Omega$, then $(\bar{x}, \bar{y}^{*}, \bar{\mu}) \in \mathcal{S}(\hyperlink{RD}{\emph{RD}_{MW}})$.
\end{itemize}
\end{thm}

\begin{proof}
Thanks to Theorem \ref{Thm3-2}, we find $y^{*} \in K^{+} \setminus \{0\}$, $u^{*} \in \partial \langle y^{*}, f \rangle(\bar{x})$, $\mu_{i} \ge 0$, $v_{i}^{*} \in \textrm{cl}^{*}\textrm{co} \Big(\medcup \Big\{\partial_{x} g_{i}(\bar{x}, v_{i}) \,\mid\, v_{i} \in \mathcal{V}_{i}(\bar{x})\Big\}\Big)$, $i = 1, 2, \dots, n$, and $b^{*} \in B_{X^{*}}$ such that
\begin{align}
  & -\Big(u^{*} + \sum_{i=1}^{n} \mu_{i} \, v^{*}_{i} + \langle y^{*}, \vartheta \rangle b^{*} \Big) \in N(\bar{x}; \Omega), \nonumber \\
  \label{4-6}
  & \mu_{i} \, \max_{v_{i} \in \mathcal{V}_{i}} g_{i}(\bar{x}, v_{i}) = 0, \quad i = 1, 2, \dots, n.
\end{align}
Putting $\bar{y}^{*} := y^{*}$ and $\bar{\mu} := (\mu_{1},\mu_{2},\dots,\mu_{n})$, we get $(\bar{y}^{*}, \bar{\mu}) \in K^{+} \setminus \{0\} \times \mathbb{R}^{n}_{+}$. Furthermore, the inclusion $v_{i} \in \mathcal{V}_{i}(\bar{x})$ means that $g_{i}(\bar{x}, v_{i}) = \max\limits_{u_{i} \in \mathcal{V}_{i}} g_{i}(\bar{x}, u_{i})$ for all $i \in \{1, 2, \dots, n\}$. Thus, it stems from (\ref{4-6}) that $\mu_{i} \, g_{i}(\bar{x}, v_{i}) = 0$, $i = 1, 2, \dots, n$. So $(\bar{x}, \bar{y}^{*}, \bar{\mu}) \in F_{MW}$.

(i) As $(f, g)$ is $\vartheta\mbox{-}\textrm{type}$ I pseudo convex on $\Omega$ at $z$ for all $z \in \Omega$, applying (i) of Theorem \ref{Thm4-1} gives us
\begin{equation*}
  \bar{f}(\bar{x}, \bar{y}^{*}, \bar{\mu}) = f(\bar{x}) \nprec \bar{f}(z, y^{*}, \mu) - \|\bar{x} - z\|\vartheta
\end{equation*}
for each $(z, y^{*}, \mu) \in F_{MW}$. Therefore $(\bar{x}, \bar{y}^{*}, \bar{\mu}) \in \mathcal{S}^{w}(\hyperlink{RD}{\textrm{RD}_{MW}})$.

(ii) As $(f, g)$ is $\vartheta\mbox{-}\textrm{type}$ II pseudo convex on $\Omega$ at $z$ for all $z \in \Omega$, applying (ii) of Theorem \ref{Thm4-1} allows us
\begin{equation*}
  \bar{f}(\bar{x}, \bar{y}^{*}, \bar{\mu}) \npreceq \bar{f}(z, y^{*}, \mu) - \|\bar{x} - z\|\vartheta
\end{equation*}
for each $(z, y^{*}, \mu) \in F_{MW}$. Therefore $(\bar{x}, \bar{y}^{*}, \bar{\mu}) \in \mathcal{S}(\hyperlink{RD}{\textrm{RD}_{MW}})$.
\end{proof}

\begin{rem}\hypertarget{Exa4-1(Rem)}{}
If in Theorem \ref{Thm4-1} and Theorem \ref{Thm4-2}, we set $\vartheta = 0$ then these theorems reduce to \cite[Theorem~4.1]{54} and \cite[Theorem~4.2]{54}.
\end{rem}

\begin{thm}\label{Thm4-3}\textbf{\textsc{(Strong Duality)}}
Let $\bar{x} \in F$ be such that the robust $\vartheta\mbox{-}\mathit{approximate}$ \emph{(KKT)} condition is satisfied at this point. Then, there exists $(\bar{y}^{*}, \bar{\mu}) \in K^{+} \setminus \{0\} \times \mathbb{R}^{n}_{+}$ such that $(\bar{x}, \bar{y}^{*}, \bar{\mu}) \in F_{MW}$. Moreover,
\begin{itemize}
  \item [\emph{(i)}] If $(f, g)$ is $\vartheta\mbox{-}\mathit{type}$ I pseudo convex on $\Omega$ at $z$ for all $z \in \Omega$, then $(\bar{x}, \bar{y}^{*}, \bar{\mu}) \in \mathcal{S}^{w}(\hyperlink{RD}{\emph{RD}_{MW}})$ and $\bar{x} \in \mathcal{S}^{w}(\hyperlink{RP}{\emph{RP}})$.
  \item [\emph{(ii)}] If $(f, g)$ is $\vartheta\mbox{-}\mathit{type}$ II pseudo convex on $\Omega$ at $z$ for all $z \in \Omega$, then $(\bar{x}, \bar{y}^{*}, \bar{\mu}) \in \mathcal{S}(\hyperlink{RD}{\emph{RD}_{MW}})$ and $\bar{x} \in \mathcal{S}(\hyperlink{RP}{\emph{RP}})$.
\end{itemize}
\end{thm}

\begin{proof}
Since $\bar{x} \in F$ satisfies the robust $\vartheta\mbox{-}\text{approximate}$ (KKT) condition, we find $y^{*} \in K^{+} \setminus \{0\}$, $u^{*} \in \partial \langle y^{*}, f \rangle(\bar{x})$, $\mu_{i} \ge 0$, $v_{i}^{*} \in \textrm{cl}^{*}\textrm{co} \Big(\medcup \Big\{\partial_{x} g_{i}(\bar{x}, v_{i}) \,\mid\, v_{i} \in \mathcal{V}_{i}(\bar{x})\Big\}\Big)$, $i = 1, 2, \dots, n$, and $b^{*} \in B_{X^{*}}$ such that
\begin{align*}
  & -\Big(u^{*} + \sum_{i=1}^{n} \mu_{i} \, v^{*}_{i} + \langle y^{*}, \vartheta \rangle b^{*} \Big) \in N(\bar{x}; \Omega), \nonumber \\
  & \mu_{i} \, \max_{v_{i} \in \mathcal{V}_{i}} g_{i}(\bar{x}, v_{i}) = 0, \quad i = 1, 2, \dots, n.
\end{align*}
Now similar to the proof of Theorem \ref{Thm4-2}, we can arrive at the result.
\end{proof}

\begin{rem}\hypertarget{Exa4-2(Rem)}{}
\begin{itemize}
  \item [(i)] Theorem \ref{Thm4-1} and Theorem \ref{Thm4-3} develop \cite[Theorem~4.18]{52} and \cite[Theorem~4.19]{52} with $Y = \mathbb{R}^{p}$.
  \item [(ii)] Theorem \ref{Thm4-1} and Theorem \ref{Thm4-3} develop \cite[Theorem~5.2]{8} and \cite[Theorem~5.3]{8} with $\Omega = X$ and $\vartheta = 0$.
\end{itemize}
\noindent Note further that our approach here is totally different from those ones presented in the aforementioned papers.
\end{rem}

We conclude this section by presenting converse duality relations between (\hyperlink{RP}{RP}) and (\hyperlink{RD}{$\textrm{RD}_{MW}$}).
\begin{thm}\label{Thm4-4}\textbf{\textsc{(Converse Duality)}}
Let $(\bar{x}, \bar{y}^{*}, \bar{\mu}) \in F_{MW}$ be such that $\bar{x} \in F$.
\begin{itemize}
   \item [\emph{(i)}] If $(f, g)$ is $\vartheta\mbox{-}\mathit{type}$ I pseudo convex on $\Omega$ at $\bar{x}$, then $\bar{x} \in \mathcal{S}^{w}(\hyperlink{RP}{\emph{RP}})$.
   \item [\emph{(ii)}] If $(f, g)$ is $\vartheta\mbox{-}\mathit{type}$ II pseudo convex on $\Omega$ at $\bar{x}$, then $\bar{x} \in \mathcal{S}(\hyperlink{RP}{\emph{RP}})$.
\end{itemize}
\end{thm}

\begin{proof}
Since $(\bar{x}, \bar{y}^{*}, \bar{\mu}) \in F_{MW}$, there exist $u^{*} \in \partial \langle \bar{y}^{*}, f \rangle(\bar{x})$, $\bar{\mu}_{i} \ge 0$, $v_{i}^{*} \in \textrm{cl}^{*}\textrm{co} \Big(\medcup \Big\{\partial_{x} g_{i}(\bar{x}, v_{i}) \,\mid\, v_{i} \in \mathcal{V}_{i}(\bar{x})\Big\}\Big)$, $i = 1, 2, \dots, n$, and $b^{*} \in B_{X^{*}}$ such that
\begin{align}\label{4-7}
  & -\Big(u^{*} + \sum_{i=1}^{n} \bar{\mu}_{i} \, v^{*}_{i} + \langle y^{*}, \vartheta \rangle b^{*} \Big) \in N(\bar{x}; \Omega),\\
  &  \bar{\mu}_{i} \, g_{i}(\bar{x}, v_{i}) \ge 0, \quad i = 1, 2, \dots, n. \nonumber
\end{align}

Let us prove (i) by contradiction. Suppose that $\bar{x} \notin \mathcal{S}^{w}(\hyperlink{RP}{\textrm{RP}})$. Therefore, there is $\hat{x} \in F$ such that $f(\hat{x}) - f(\bar{x}) + \|\hat{x} - \bar{x}\| \vartheta \in - \textrm{int}\hspace{.4mm}K$. The latter inclusion provides $\langle \bar{y}^{*}, f(\hat{x}) - f(\bar{x}) + \|\hat{x} - \bar{x}\| \vartheta \rangle < 0$. By the $\vartheta\mbox{-}\textrm{type}$ I pseudo convex on $\Omega$ at $\bar{x}$, we infer from this inequality that there exists $w \in -N(\bar{x}; \Omega)^{+}$ such that
\begin{align*}
&\langle u^{*}, w \rangle +  \langle y^{*}, \vartheta \rangle \|\hat{x} - \bar{x}\| < 0, \\
&\|w\| \le \|\hat{x} - \bar{x}\|.
\end{align*}

Moreover, from (\ref{4-7}) we have for $w$
\begin{equation*}
  \langle u^{*}, w \rangle + \sum_{i=1}^{n} \bar{\mu}_{i} \, \langle v^{*}_{i}, w \rangle + \langle y^{*}, \vartheta \rangle \langle b^{*}, w \rangle \ge 0.
\end{equation*}
So, the above relationships entail that
\begin{equation*}
  \sum_{i=1}^{n} \bar{\mu}_{i} \, \langle v^{*}_{i}, w \rangle > 0.
\end{equation*}
Now argue as in Theorem \ref{Thm3-4}(i)'s proof, one can arrive at the result.

The proof of (ii) is similar to that of (i), so we omit the corresponding details.
\end{proof}

\begin{rem}\hypertarget{Exa4-3(Rem)}{}
If in Theorem \ref{Thm4-4}, we set $\vartheta = 0$ then this theorem reduces to \cite[Theorem~4.4]{54}.
\end{rem}

\pdfbookmark[section]{\refname}{bibliography}
\printbibliography
\small

\end{document}